\documentclass[preprint,12pt]{elsarticle}

\usepackage{amssymb}
\usepackage{amsthm}
\usepackage{graphics}
\usepackage{epsfig}
\usepackage{graphicx}
\usepackage{color}
\usepackage{hyperref}
\usepackage{a4wide}
\usepackage{amsmath}
\usepackage{amsfonts}
\usepackage{enumerate}
\newcommand{\pf}{\noindent {\bf Proof: }}

\newtheorem*{theoremaux}{Theorem \theoremauxnum}
\gdef\theoremauxnum{1}

\newtheorem{lemma}{\bf Lemma}[section]
 
\newtheorem{theorem}{\bf Theorem}[section]
\newtheorem{proposition}[lemma]{\bf Proposition}
\newtheorem{corollary}[lemma]{\bf Corollary}
\newtheorem{definition}{\bf Definition}[section]
\newtheorem{remark}{\bf Remark}[section]

\journal{~}

\begin{document}

\begin{frontmatter}



\title{Classification of Cayley Rose Window Graphs}



\author{Angsuman Das\corref{cor1}}
\ead{angsuman.maths@presiuniv.ac.in}

\author{Arnab Mandal}
\ead{arnab.maths@presiuniv.ac.in}

\address{Department of Mathematics,\\ Presidency University, Kolkata\\ 86/1, College Street, Kolkata - 7000073\\West Bengal, India}
\cortext[cor1]{Corresponding author}


\begin{abstract}
	Rose window graphs are a family of tetravalent graphs, introduced by Steve Wilson. Following it, Kovacs, Kutnar and Marusic classified the edge-transitive rose window graphs and Dobson, Kovacs and Miklavic characterized the vertex transitive rose window graphs.	In this paper, we classify the Cayley rose window graphs.
\end{abstract}

\begin{keyword}
	{vertex-transitive \sep  regular subgroup \sep rose window graph.}
	\MSC[2010]  {\bf 05C75, 05E18}
	
\end{keyword}

\end{frontmatter}


\section{Introduction}
{\it Rose window graphs} were introduced in \cite{wilson} in the following way:
{\definition Given natural numbers $n\geq 3$ and $1\leq a,r\leq n-1$, the Rose Window graph $R_n(a,r)$ is defined to be the graph with vertex set $V=\{A_i,B_i:i \in \mathbb{Z}_n\}$ and four kind of edges: $A_iA_{i+1}$ $(\mathbf{rim}~ edges)$, $A_iB_i$ $(\mathbf{ inspoke}~ edges)$, $A_{i+a}B_i$ $(\mathbf{ outspoke}~ edges)$ and $B_iB_{i+r}$ $(\mathbf{ hub}~ edges)$, where the addition of indices are done modulo $n$.}
 
In the introductory paper \cite{wilson}, author's initial interest in rose window graphs arose in the context of graph embeddings into surfaces. The author conjectured that rose window graphs are edge-transitive if and only if it belongs to the one of the four families given in Theorem \ref{ET-theorem}. The conjecture was proved by Kovacs {\it et. al.} in \cite{kovacs}. In particular, they proved that 

{\theorem \cite{kovacs} \label{ET-theorem} A rose window graph is edge-transitive if and only if it belongs to one of the four families:
	\begin{enumerate}
		\item $R_n(2,1)$.
		\item $R_{2m}(m\pm 2,m\pm 1)$
		\item $R_{12m}(\pm (3m+2),\pm (3m-1))$ and $R_{12m}(\pm (3m-2),\pm (3m+1))$.
		\item $R_{2m}(2b,r)$, where $b^2\equiv \pm 1$(mod $m$), $2\leq 2b\leq m$, and $r \in \{1,m-1\}$ is odd.
	\end{enumerate}
	} 
A similar characterization for vertex-transitive graphs was proved in \cite{vt-rw}:
{\theorem\cite{vt-rw} \label{VT-theorem} A rose window graph $R_n(a,r)$ is vertex-transitive if and only if it belongs to one of the following families:
	\begin{enumerate}
		\item $R_n(a,r)$, where $r^2\equiv \pm 1(mod~ n)$ and $ra\equiv \pm a(mod~ n)$.
		\item $R_{4m}(2m,r)$, where $r$ is odd and $(r^2+2m)\equiv \pm 1(mod~ 4m)$.
		\item $R_{2m}(m\pm 2,m\pm 1)$
		\item $R_{12m}(\pm (3m+2),\pm (3m-1))$ and $R_{12m}(\pm (3m-2),\pm (3m+1))$.
		\item $R_{2m}(2b,r)$, where $b^2\equiv \pm 1$(mod $m$), $2\leq 2b\leq m$, and $r \in \{1,m-1\}$ is odd.
	\end{enumerate}    }

As a Cayley graph is always vertex-transitive, a natural question to ask is to characterize the rose-window graphs which are also Cayley graphs. For that, it is sufficient to look for Cayley graphs only in the $5$ families mentioned in Theorem \ref{VT-theorem}. The main goal of this paper is finding an answer to this question. In particular, we prove the following theorem:
\begin{theorem}\label{cayley-iff} A rose-window graph $R_n(a,r)$ is Cayley if and only if one of the following holds:
\begin{enumerate}
	\item $R_n(a,r)$, where $r^2\equiv \pm 1(mod~ n)$ and $ra\equiv \pm a(mod~ n)$.
	\item $R_{4m}(2m,r)$, where $r$ is odd and $\mathbf{(r^2+2m)\equiv 1(mod~ 4m)}$.
	\item $R_{2m}(m\pm 2,m\pm 1)$ where $\mathbf{m}$ {\bf is a multiple of $\mathbf{2}$ or $\mathbf{3}$}.
	\item $R_{12m}(\pm (3m+2),\pm (3m-1))$ and $R_{12m}(\pm (3m-2),\pm (3m+1))$ where $\mathbf{m\not\equiv 0(mod ~4)}$.
	\item $R_{2m}(2b,r)$, where $b^2\equiv \pm 1$(mod $m$), $2\leq 2b\leq m$, and $r \in \{1,m-1\}$ is odd.\qed
\end{enumerate}	
\end{theorem}	

Before stating the proof, we note a few generic automorphisms and other properties of $R_n(a,r)$. Other automorphisms, specific to any particular family of rose window graphs, will be introduced whenever they are needed.
\begin{enumerate}
	\item Define $\tau: V \rightarrow V$ by $\tau(A_i)=A_{-i}$ and $\tau(B_i)=B_{-i}$. Clearly $\tau$ is an automorphism with $\tau^2=\mathsf{id}$ and hence $R_n(a,r)\cong R_n(-a,r)$. 
	\item $R_n(a,r)= R_n(a,-r)$.
	\item Define $\rho: V \rightarrow V$ by $\rho(A_i)=A_{i+1}$ and $\rho(B_i)=B_{i+1}$; and $\mu: V \rightarrow V$ by $\mu(A_i)=A_{-i}$ and $\mu(B_i)=B_{-a-i}$. Clearly $\rho$ and $\mu$ are automorphisms. As $\rho^n=\mu^2=\mathsf{id}$ and $\mu\rho\mu=\rho^{-1}$, we have $\langle\rho,\mu \rangle\cong D_n$.
	\item If $(n,r)=1$, then $\zeta : V\rightarrow V$ given by $\zeta(A_i)=B_{-ir^{-1}}$ and $\zeta(B_i)=A_{-ir^{-1}}$ is an automorphism and hence $R_n(a,r)\cong R_n(ar^{-1},r^{-1})$.
\end{enumerate}

{\remark \label{r,a<=n/2-remark} In view of the first two observations, it is enough to study $R_n(a,r)$ for $1\leq a,r \leq \lceil\frac{n}{2} \rceil$.} 

The main theorem, which is repeatedly used in the proofs throughout the paper, is the following:

{\proposition A vertex-transitive graph $G$ is Cayley if and only if $\mathsf{Aut}(G)$ has a subgroup $H$ which acts regularly on the vertices of $G$. In particular, non-identity elements of $H$ do not stabilize any vertex.\qed}

{\remark \label{regular-transitive} In this context, it is to be noted that if a group of order $n$ acts transitively on a set of order $n$, then the action is regular.}

\section{Family-1 [$R_n(a,r)$: $r^2\equiv \pm 1(mod~ n)$ and $ra\equiv \pm a(mod~ n)$]}
If $r^2\equiv \pm 1(mod~ n)$ and $ra\equiv \pm a(mod~ n)$, then $\delta : V\rightarrow V$ given by $\delta(A_{i})=B_{ri}$ and $\delta(B_i)=A_{ri}$ is an automorphism. For proof, see Lemma 2 \cite{wilson} or Lemma 3.7 \cite{vt-rw}. If $r^2\equiv 1(mod~ n)$, then $\delta^2=\mathsf{id}$ and if $r^2\equiv -1(mod~ n)$, then $\delta^2=\tau$, i.e., $\delta$ is of order $4$.
{\theorem \label{r^2=1} If $r^2\equiv 1(mod~ n)$ and $ra\equiv \pm a(mod~ n)$, then $R_n(a,r)$ is a Cayley graph.}\\
\pf Since $R_n(a,r)=R_n(a,-r)$, without loss of generality, we can assume that $ra\equiv - a(mod~ n)$. Consider $\rho$ and $\delta$ as defined above. We have $\rho^n=\delta^2=\mathsf{id}$ and $\delta\rho\delta=\rho^r$. Define $$H=\langle\rho,\delta\rangle = \langle\rho,\delta: \rho^n=\delta^2=\mathsf{id};\delta\rho\delta=\rho^r\rangle$$
$$~~~~~~~~~~~~~~~~~~~~~~~~~~~=\{\mathsf{id},\rho,\rho^2,\ldots,\rho^{n-1},\delta,\rho\delta,\rho^2\delta,\ldots,\rho^{n-1}\delta \}.$$
Clearly, $H$ is a subgroup of $\mathsf{Aut}(R_n(a,r))$. It suffices to show that $H$ acts regularly on $R_n(a,r)$. For that we observe that \begin{itemize}
	\item $\rho^j(A_i)=A_{i+j}$ and $\rho^j(B_i)=B_{i+j}$, and
	\item $\rho^j\delta(A_i)=B_{ri+j}$ and $\rho^j\delta(B_i)=A_{ri+j}$.
\end{itemize}  As $gcd(r,n)=1$, the map $i\mapsto ri+j$ is a bijection on $\{0,1,\ldots,n-1\}$. Thus $H$ acts transitively on $R_n(a,r)$. It is also clear from the construction of $H$, that for any pair of vertices in $R_n(a,r)$, there exists a unique element in $H$ which maps one to the other. Hence, $R_n(a,r)$ is a Cayley graph.\qed

{\lemma If $r^2\equiv -1(mod~ n)$ and $ra\equiv \pm a(mod~ n)$, then $n$ is even, $a$ is odd and $n=2a$. }\\
\pf Let $p$ be an odd prime factor of $n$ such that $p^i|n$ and $p^{i+1}\nmid n$. Then $r^2\equiv -1(mod ~p^i)$ and $r^2\equiv -1(mod ~p)$. Again,  $p^i|a(r\pm 1)$, i.e., $p|a(r\pm 1)$. If $p|(r\pm 1)$, then $r^2\equiv 1 (mod~p)$, a contradiction, as $-1\not\equiv 1(mod~p)$. Thus for all odd prime factors $p$ of $n$, we have $p^i|a$. Hence, if $n$ is odd, then $n=a$, a contradiction (See Remark \ref{r,a<=n/2-remark}). Thus $n$ is even.

We claim that $2|n$ but $4\nmid n$. Because if $4|n$, then $r^2\equiv -1(mod~4)$. However, there does not exist any such $r$. Thus $n$ is $2$ times the product of some odd primes. Also, all the odd prime factors of $n$ are also factors of $a$, as seen above. Thus, if $2|a$, then $n=a$, a contradiction (See Remark \ref{r,a<=n/2-remark}). Thus $2\nmid a$ and hence $a$ is odd and $n=2a$.\qed

{\theorem If $r^2\equiv -1(mod~ n)$ and $ra\equiv \pm a(mod~ n)$, then $R_n(a,r)$ is a Cayley graph.}\\
\pf Let $\alpha=\rho^2; \beta=\rho\delta^2; \gamma=\mu\delta$. Clearly, $\alpha,\beta,\gamma \in \mathsf{Aut}(R_n(a,r))$. It can be easily checked that $\beta\alpha=\alpha^{-1}\beta;\gamma\alpha=\alpha^{-r}\gamma$ and $\gamma^2=\alpha^{\frac{a-1}{2}}\beta$. Define
 $$H=\langle\alpha,\beta,\gamma :\alpha^{n/2}=\beta^2=\gamma^4=\mathsf{id}; \beta\alpha=\alpha^{-1}\beta;\gamma\alpha=\alpha^{-r}\gamma;\gamma^2=\alpha^{\frac{a-1}{2}}\beta  \rangle$$
 $$=\{\alpha^i\beta^j\gamma^k : 0\leq i < n/2, 0\leq j,k\leq 1  \}~~~~~~~~~~~~~~~~~~~~~~~~~~~~~~~~~~~~~$$ 
\noindent Note that, from the above lemma, $n/2$ and $(a-1)/2$ are positive integers. We claim that the elements in $H$ are distinct. If not, suppose $$\alpha^{i_1}\beta^{j_1}\gamma^{k_1}=\alpha^{i_2}\beta^{j_2}\gamma^{k_2}, \mbox{ where } 0\leq i_1,i_2 < n/2, 0\leq j_1,j_2\leq 1,0\leq k_1,k_2\leq 1,$$ i.e., $$\beta^{-j_2}\alpha^{i_1-i_2}\beta^{j_1}=\gamma^{k_2-k_1}, \mbox{ where }k_2-k_1=0 \mbox{ or }1.$$ Now, as $\gamma=\mu\delta$ flips $A_i$'s and $B_j$'s, and $\alpha,\beta$ maps $A_i$'s to $A_j$'s and $B_i$'s to $B_j$'s, $k_2-k_1$ must be $0$, i.e., $k_1=k_2$. Thus, we have $$\alpha^{i_1-i_2}=\beta^{j_2-j_1}, \mbox{ where } j_2- j_1=0 \mbox{ or }1.$$  If $j_2-j_1=1$, then $\alpha^{i_1-i_2}=\beta=\rho\delta^2$. But $\alpha^{i_1-i_2}(A_0)=A_{2(i_1-i_2)}$ (even index) and $\rho\delta^2(A_0)=A_1$ (odd index). Hence, $j_2-j_1=0$, i.e., $j_1=j_2$. This implies $\alpha^{i_1-i_2}=\mathsf{id}$ and as a result $i_1=i_2$. Thus the elements of $H$ are distinct and $|H|=n/2\times 2 \times 2=2n$.

\noindent We claim that $H$ acts transitively on $R_n(a,r)$. It suffices to show that the stabilizer of $A_0$ in $H$, $\mathsf{Stab}_H(A_0)=\{\mathsf{id} \}$.

Let $\alpha^i\beta^j\gamma^k \in \mathsf{Stab}_H(A_0)$, i.e., $\alpha^i\beta^j\gamma^k(A_0)=A_0$. Since, $\gamma$ flips $A_i$'s and $B_j$'s, and $\alpha,\beta$ do not, we have $k=0$. Thus, $\alpha^i\beta^j(A_0)=A_0$. If $j=1$, then $\alpha^i\beta(A_0)=\alpha^i\rho\delta^2(A_0)=\rho^{1+2i}\delta^2(A_0)=A_0$, i.e., $A_{1+2i}=A_0$, a contradiction, as the parity of indices on both sides does not match. Thus, $j=0$ and we have $\alpha^i(A_0)=A_0$. But this implies $A_{2i}=A_0$, i.e., $i=0$. Hence $\mathsf{Stab}_H(A_0)=\{\mathsf{id} \}$.

Finally, in view of Remark \ref{regular-transitive}, $H$ acts regularly on $R_n(a,r)$ and hence $R_n(a,r)$ is a Cayley graph.\qed


\section{Family-2 [$R_{4m}(2m,r)$: $r$ is odd and $(r^2+2m)\equiv \pm 1(mod~ 4m)$]}
{\proposition \label{family-2-gamma-definition} If $n$ is divisible by $4$, $r$ is odd, $a=n/2$ and $(r^2+n/2)\equiv \pm 1(mod~ n)$, then
\begin{itemize}
	\item $gcd(r,n)=1$.
		\item If $\gamma: V \rightarrow V$ be defined by $\gamma(A_i)=B_{ri}$ and $\gamma(B_i)=A_{(r+a)i}$, then $\gamma \in \mathsf{Aut}(R_n(a,r))$.
\end{itemize}			
	}
	\pf Let $n=4m$ and $a=2m$, and let if possible, $gcd(r,n)=l>1$. As $r$ is odd, $l|m$. Thus $r=lt$ and $m=ls$ for some $s,t \in \mathbb{N}$. Thus $n=4ls, a=2ls$ and $r=lt$. Now $(r^2+n/2)\equiv\pm 1(mod~ n)$ implies $l^2t^2+2ls\equiv \pm 1(mod~ 4ls)$, which in turn implies $l|(l^2t^2+2ls\pm 1)$, i.e., $l|1$, a contradiction. Thus $gcd(r,n)=1$.
	
	$\gamma$, as defined above, has been shown to be in $\mathsf{Aut}(R_n(a,r))$ in Lemma 3.8 \cite{vt-rw}. 	\qed
{\proposition \label{family-2-zeta-definition} If $n$ is divisible by $4$, $r$ is odd, $a=n/2$ and $(r^2+n/2)\equiv 1(mod~ n)$, then
\begin{itemize}
	\item $r^{-1}=r+a$ (mod $n$)
	\item $\zeta \in \mathsf{Aut}(R_n(a,r))$ (defined before) takes the following form:  $\zeta(A_i)=B_{-(r+a)i}$ and $\zeta(B_i)=A_{-(r+a)i}$, and $\zeta^4=\mathsf{id}$.
\end{itemize}	}

	\pf $r(r+a)\equiv r^2+ar\equiv 1-a+ar\equiv 1+a(r-1)\equiv 1 ~(mod ~n)$. The last equivalence holds as $r$ is odd and $a=n/2$. Thus $r^{-1}=r+a$ (mod $n$). The form of $\zeta$ follows immediately from the fact that $r^{-1}=r+a$ (mod $n$).\qed

{\theorem \label{family-2-cayley-theorem} If $n$ is divisible by $4$, $r$ is odd, $a=n/2$ and $(r^2+n/2)\equiv 1(mod~ n)$, then $R_n(a,r)$ is a Cayley graph.}\\
\pf Let $\alpha=\rho^2,\beta=\rho\mu$ and $\sigma=\gamma \zeta^2$, where $\gamma$ and $\zeta$ are as defined in Propositions \ref{family-2-gamma-definition} and \ref{family-2-zeta-definition}. It can be easily checked that $\sigma(A_i)=B_{(r+a)i}$ and $\sigma(B_i)=A_{ri}$; $\alpha^{n/2}=\beta^2=\sigma^2=\mathsf{id}$; $\beta\alpha\beta=\alpha^{-1},\sigma\alpha\sigma=\alpha^r, (\beta\sigma)^2=\alpha^{\frac{a-r+1}{2}}$. Define
$$H=\langle\alpha,\beta,\sigma :\alpha^{n/2}=\beta^2=\sigma^2=\mathsf{id}; \beta\alpha\beta=\alpha^{-1},\sigma\alpha\sigma=\alpha^r, (\beta\sigma)^2=\alpha^{\frac{a-r+1}{2}}  \rangle$$
$$=\{\alpha^i\beta^j\sigma^k : 0\leq i < n/2, 0\leq j,k\leq 1 \}~~~~~~~~~~~~~~~~~~~~~~~~~~~~~~~~~~~~~~~~$$ 
\noindent We claim that the elements in $H$ are distinct. If not, suppose $$\alpha^{i_1}\beta^{j_1}\sigma^{k_1}=\alpha^{i_2}\beta^{j_2}\sigma^{k_2}, \mbox{ where } 0\leq i_1,i_2 < n/2, 0\leq j_1,j_2,k_1,k_2\leq 1,$$ i.e., $$\alpha^{i_1-i_2}\beta^{j_1}\sigma^{k_1-k_2}=\beta^{j_2}, \mbox{ where }k_1-k_2=0 \mbox{ or }1.$$ Now, as $\sigma$ flips $A_i$'s and $B_j$'s, and $\alpha,\beta$ maps $A_i$'s to $A_j$'s and $B_i$'s to $B_j$'s, $k_1-k_2$ must be $0$, i.e., $k_1=k_2$. Thus, we have $$\alpha^{i_1-i_2}=\beta^{j_2-j_1}, \mbox{ where } j_2- j_1=0 \mbox{ or }1.$$  Since, $\alpha$ maintains the parity of indices and $\beta$ flips the parity of indices of $A_i$'s and $B_i$'s, $j_2-j_1$ is even, i.e., $j_1=j_2$. This implies $\alpha^{i_1-i_2}=\mathsf{id}$ and as a result $i_1=i_2$. Thus the elements of $H$ are distinct and $|H|=n/2\times 2 \times 2=2n$.

\noindent We claim that $H$ acts transitively on $R_n(a,r)$. In order to prove it, we show that the orbit of $A_0$, $\mathcal{O}_{A_0}$, under the action of $H$ is the vertex
set of $R_n(a, r)$. By orbit-stabilizer theorem, we get
$$|\mathcal{O}_{A_0}|=\dfrac{|H|}{|\mathsf{Stab}_H (A_0)|}.$$
As the number of vertices in $R_n(a,r)$ is $2n$ and $|H| = 2n$, it is enough to
show that $\mathsf{Stab}_H (A_0) = \{\mathsf{id}\}$. Let $\alpha^i\beta^j\sigma^k$ be an arbitrary element of $H$ which stabilizes $A_0$, i.e., $\alpha^i\beta^j\sigma^k(A_0)=A_0$, with $0\leq i < n/2, 0\leq j,k\leq 1$. Now, as $\sigma$ flips $A_i$'s and $B_j$'s, and $\alpha,\beta$ maps $A_i$'s to $A_j$'s and $B_i$'s to $B_j$'s, $k=0$. Thus $\alpha^i\beta^j(A_0)=A_0$, i.e., $\alpha^{-i}(A_0)=\beta^j(A_0)$. Since, $\alpha$ maintains the parity of indices and $\beta$ flips the parity of indices of $A_i$'s and $B_i$'s, $j=0$ and hence $i=0$. Thus $\mathsf{Stab}_H (A_0) = \{\mathsf{id}\}$.

Finally, in view of Remark \ref{regular-transitive}, $H$ acts regularly on $R_{n}(a,r)$ and hence $R_{n}(a,r)$ is a Cayley graph.\qed

In {\bf Family 2}, if $(r^2+n/2)\equiv -1(mod~ n)$, we will show that $R_n(a,r)$ is not a Cayley graph. In order to prove it, we recall a few observations and results.
{\remark \label{edge-orbit-remark} It was noted in \cite{wilson} and \cite{vt-rw}, that $R_n(a,r)$ has either one or two or three edge orbits. If it has one edge orbit,  then by definition, it is edge transitive, as in Theorem \ref{ET-theorem}. If $R_n(a, r)$ has two edge orbits, then one orbit consists of rim and hub edges, and the other consists of spoke edges. If $R_n (a, r)$ has three orbits on edges, then the first one consists of rim edges, the second one consists of hub edges, and the third one consists of spoke edges.}
	
 As {\bf Family} $\mathbf{3,4,5}$ in Theorem \ref{VT-theorem} are also edge transitive, they have only one edge orbit. On the other hand, family $1$ and $2$ in Theorem \ref{VT-theorem}, have two edge orbits, as evident from Remark \ref{edge-orbit-remark} and Theorem \ref{rim-to-hub-edge-theorem}.

{\theorem[Theorem 2.3,\cite{vt-rw}] \label{rim-to-hub-edge-theorem} There is an automorphism	of $R_n(a,r)$ sending every rim edge to a hub edge and vice-versa if and only if one of the following holds:
	\begin{enumerate}
		\item $a\neq n/2$, $r^2\equiv 1 (mod ~n)$ and $ra\equiv \pm a (mod ~n)$;
		\item $a= n/2$, $r^2\equiv \pm 1 (mod ~n)$ and $ra\equiv \pm a (mod ~n)$;
		\item $n$ is divisible by $4$, $gcd(n, r) = 1$, $a = n/2$ and $(r^2 + n/2)\equiv \pm 1 (mod ~n)$.
	\end{enumerate}		}
	
{\corollary[Corollary 3.9,\cite{vt-rw}] \label{family-2-full-aut-group} If $n$ is divisible by $4$, $r$ is odd, $a = n/2$ and $(r^2 + n/2)\equiv \pm 1 (mod ~n)$, then the automorphism group of $R_n (a, r)$ has two edge orbits and the full automorphism group of $R_n (a, r)$, $\mathsf{Aut}(R_n(a,r))=\langle \rho,\mu,\gamma \rangle$, where $\gamma$ is as defined in Proposition \ref{family-2-gamma-definition}.}

{\theorem If $n$ is divisible by $4$, $r$ is odd, $a=n/2$ and $(r^2+n/2)\equiv -1(mod~ n)$, then $R_n(a,r)$ is not a Cayley graph.}\\
\pf As evident from Corollary \ref{family-2-full-aut-group}, the full automorphism group of $R_n (a, r)$ is given by
$$\mathsf{Aut}(R_n(a,r))=\langle \rho,\mu,\gamma: \rho^n=\mu^2=\gamma^4=\mathsf{id}; \mu\rho\mu=\rho^{-1},\gamma\mu=\rho^a\mu\gamma,\gamma\rho=\rho^{r-a}\mu\gamma^3 \rangle.$$
One can easily check the relations between the generators starting from the definition and conclude that $|\mathsf{Aut}(R_n(a,r))|=n\times 2 \times 4=8n$. If possible, let $R_n(a,r)$ be a Cayley graph with a regular subgroup $H$ of $\mathsf{Aut}(R_n(a,r))$ and $|H|=2n$.

Let $K=\langle \gamma \rangle$. Then $|K|=4$ and $H \cap  K$ is a subgroup of $K$. As $\gamma^2(A_0)=A_0$, i.e., $\gamma^2$ has a fixed point, $\gamma^2\not\in H$. Thus $H\cap K=\{\mathsf{id}\}$ and $$|HK|=\dfrac{|H||K|}{|H\cap K|}=8n.$$
Hence $\mu \in \mathsf{Aut}(R_n(a,r))=HK$. Thus $\mu =hk$, where $h \in H$ and $k \in K=\{\mathsf{id},\gamma,\gamma^2,\gamma^3\}$. If $k=\mathsf{id}$, then $\mu=h\in H$. But as $\mu(A_0)=A_0$, i.e., $\mu$ has a fixed point, $\mu \not\in H$. Thus $k\neq \mathsf{id}$.\\
If $k=\gamma^2$, then $\mu=h\gamma^2$, i.e., $h=\mu\gamma^2 \in H$. But as $\mu\gamma^2(A_0)=A_0$, $\mu\gamma^2 \not\in H$ and hence $k\neq \gamma^2$.\\
If $k=\gamma$, then $\mu\gamma^{-1}=h$, i.e., $h^{-2}=(\gamma\mu)^2=\rho^a\gamma^2 \in H$. But, as $\rho^a\gamma^2(A_{a/2})=A_{a/2}$, by similar argument, $k\neq \gamma$.\\
If $k=\gamma^3$, then $h^2=(\mu\gamma)^2=\rho^a\gamma^2 \in H$. By similar argument as above, $k\neq \gamma^3$.\\
As all the four possible choices of $k \in K$ leads to contradiction, we conclude that there does not exist any regular subgroup $H$ of $\mathsf{Aut}(R_n(a,r))$ and hence $R_n(a,r)$ is not a Cayley graph. \qed 

\section{Family-3 [$R_{2m}(m\pm 2,m\pm 1)$]}
As $m+2\equiv -(m-2) ~(mod ~2m)$ and $m+1\equiv -(m-1)~ (mod ~2m)$, it suffices to check the family $R_{2m}(m-2,m-1)$. It was proved in Section 3.2 of \cite{rotary}, that $$G:=\mathsf{Aut}(R_{2m}(m-2,m-1))=\langle \rho,\mu, \varepsilon_0,\varepsilon_1,\ldots,\varepsilon_{m-1} \rangle=K\rtimes \langle \rho\varepsilon_0,\mu\rho^m \rangle \cong \mathbb{Z}_2^m \rtimes D_m,$$ where $K=\langle \varepsilon_0,\varepsilon_1,\ldots,\varepsilon_{m-1} \rangle \cong \mathbb{Z}_2^m$, $D_m$ is the dihedral group and $\varepsilon_i$ is the involution given by $(A_i,B_{i-1})(A_{i+m},B_{i-1+m})(A_{i+1},B_{i+m})(A_{i+1+m},B_{i})$. Thus $|G|=2^{m+1}m$. One can easily check that the following relations between the generators hold: $$\varepsilon_i\varepsilon_j=\varepsilon_j\varepsilon_i;~ \varepsilon_i\rho^m=\rho^m\varepsilon_i;~\mu\varepsilon_i=\varepsilon_{m-1-i}\mu;~$$
$$\rho\varepsilon_i=\varepsilon_{i+1}\rho, \forall i,j \in \{0,1,\ldots,m-1\} \mbox{ and }\varepsilon_0\varepsilon_1\cdots\varepsilon_{m-1}=\rho^m$$ where the addition of indices of $\varepsilon_i$'s are done modulo $m$. Using this relations, it is easy to see that $\circ (\rho\varepsilon_i)=m$ and $\circ (\mu\rho^i)=2$. 

It follows from definition that $\rho^{2i}\mu,\varepsilon_0,\varepsilon_1,\ldots,\varepsilon_{i-2},\varepsilon_{i+1},\ldots,\varepsilon_{m-1} \in \mathsf{Stab}_G(A_i)$. Again, using the relations between generators, we get $|\langle \rho^{2i}\mu,\varepsilon_0,\varepsilon_1,\ldots,\varepsilon_{i-2},\varepsilon_{i+1},\ldots,\varepsilon_{m-1}\rangle|=2^{m-1}$. Now, as $R_{2m}(m-2,m-1)$ is a vertex transitive graph, by orbit-stabilizer theorem, it follows that $|G|/|\mathsf{Stab}_G(A_i)|=2\times 2m,$ i.e., $|\mathsf{Stab}_G(A_i)|=\frac{2^{m+1}m}{4m}=2^{m-1}$. Thus, we have $$\mathsf{Stab}_G(A_i)=\langle \rho^{2i}\mu,\varepsilon_0,\varepsilon_1,\ldots,\varepsilon_{i-2},\varepsilon_{i+1},\ldots,\varepsilon_{m-1} \rangle.$$ Similarly, it follows that $$\mathsf{Stab}_G(B_i)=\langle \rho^{m-2+2i}\mu,\varepsilon_0,\varepsilon_1,\ldots,\varepsilon_{i-1},\varepsilon_{i+2},\ldots,\varepsilon_{m-1} \rangle.$$

{\theorem $R_{2m}(m-2,m-1)$ is a Cayley graph, if $m$ is even.}\\
\pf In this case, $n=2m$, $a=m-2$ and $r=m-1$. Now, if $m$ is even, we have $$r^2=(m-1)^2=m^2-2m+1\equiv 1 (mod~2m)\equiv 1(mod~n) \mbox{ and}$$ $$ra=(m-1)(m-2)=m^2-3m+2\equiv -m+2(mod~2m)\equiv -a(mod~n).$$ Thus, if $m$ is even, $R_{2m}(m-2,m-1)$ is a subfamily of {\bf Family-1} and as a result, $R_{2m}(m-2,m-1)$ is a Cayley graph.\qed

{\theorem $R_{2m}(m-2,m-1)$ is a Cayley graph, if $m$ is an odd multiple of $3$.}\\
\pf Let $m=3l$. For $i=0,1,2$, denote by $\Sigma_i$, the product of all $\varepsilon_j$'s such that $j\neq i~(mod~3)$. Note that $\Sigma_i\Sigma_j=\Sigma_k$ for distinct $i,j,k$'s in $\{0,1,2\}$ and $\circ(\Sigma_i)=2$.

 Let $\alpha=\rho^2,\beta=\Sigma_0$ and $\gamma=\Sigma_1$. It can be easily checked that $\beta\alpha=\alpha\gamma,\gamma\alpha=\alpha\beta\gamma$ and $\beta\gamma=\gamma\beta$. Define $$H=\langle\alpha,\beta,\gamma:\circ(\alpha)=m,\circ(\beta)=\circ(\gamma)=2;\beta\alpha=\alpha\gamma,\gamma\alpha=\alpha\beta\gamma,\beta\gamma=\gamma\beta \rangle.$$ Thus, any element of $H$ can be expressed as $\alpha^i\beta^j\gamma^k$ where $0\leq i\leq m-1,0\leq j,k\leq 1$, i.e., $|H|\leq 4m$.\\ 
 {\it Claim 1:} $|H|=4m$. \\
 {\it Proof of Claim 1:} If not, there exist $0\leq i_1,i_2\leq m-1,0\leq j_1,j_2,k_1,k_2\leq 1$ such that $\alpha^{i_1}\beta^{j_1}\gamma^{k_1}=\alpha^{i_2}\beta^{j_2}\gamma^{k_2}$, i.e., $$\rho^{2(i_1-i_2)}=\alpha^{i_1-i_2}=\beta^{j_2-j_1}\gamma^{k_2-k_1} ~(\mbox{as }\beta\gamma=\gamma\beta).$$ If $j_2-j_1=k_2-k_1=0$, then $i_1=i_2$ (since, $\circ(\rho)=2m$) and as a result the claim is true. However, if any one or both of $j_2-j_1$ or $k_2-k_1$ is $1$, then the right hand side is an element of order $2$. As a result, the left hand side must be an element of order $2$, which implies $2(i_1-i_2)=m$. However, as $m$ is odd, this can not hold. As a result, the claim is true, i.e., $|H|=4m$.
 
 As in proof of Theorem \ref{family-2-cayley-theorem}, it is enough to
 show that $\mathsf{Stab}_H (A_0) = \{\mathsf{id}\}$. Let $\alpha^i\beta^j\gamma^k\in \mathsf{Stab}_H (A_0)$, i.e., $\alpha^i\beta^j\gamma^k(A_0)=A_0$ for some $i,j,k$ with $0\leq i\leq m-1$, $0\leq j,k\leq 1$. Therefore, 
 \begin{equation}\label{family-3-odd-multiple-of-3-equation}
  \beta^j\gamma^k(A_0)=A_{2m-2i}
 \end{equation}
 {\it Claim 2:} $k=0$.\\
 {\it Proof of Claim 2:} If not, let $k=1$, i.e., $\beta^j\gamma(A_0)=A_{2m-2i}$. Note that \begin{itemize}
 	\item both $\varepsilon_0$ and $\varepsilon_{m-1}$ occurs in the expression of $\gamma$, and
 	\item all $\varepsilon_i$'s except $\varepsilon_0$ and $\varepsilon_{m-1}$ stabilizes $A_0$.
 \end{itemize}
 Thus $A_{2m-2i}=\beta^j\gamma(A_0)=\beta^j\varepsilon_{m-1}\varepsilon_0(A_0)=\beta^j\varepsilon_{m-1}(B_{2m-1})=\beta^j(A_m).$
 If $j=0$, then we have $A_m=A_{2m-2i}$, which is a contradiction, due to mismatch of parity of indices. If $j=1$, then we have $\beta(A_m)=A_{2m-2i}$. Note that 
 \begin{itemize}
 	\item $\mathsf{Stab}_G (A_0) =\mathsf{Stab}_G (A_m) =\langle \mu,\varepsilon_1,\varepsilon_2,\ldots,\varepsilon_{m-2} \rangle$. 
 	\item $\varepsilon_0$ does not occur in the expression of $\beta$, but $\varepsilon_{m-1}$ occur in the expression of $\beta$.
 \end{itemize}
 Thus, we have $A_{2m-2i}=\beta(A_m)=\varepsilon_{m-1}(A_m)=B_{2m-1}$, a contradiction. Hence for $k=1$, both $j=0$ or $j=1$ leads to a contradiction, and as a result $k=0$.
 
 Thus, from Equation \ref{family-3-odd-multiple-of-3-equation}, we have $\beta^j(A_0)=A_{2m-2i}$. If $j=1$, then $A_{2m-2i}=\beta(A_0)=\varepsilon_{m-1}(A_0)=B_{m-1}$, a contradiction. Thus, $j=0$ and hence we have $A_0=A_{2m-2i}$ i.e., $2m\equiv 2i~(mod~2m)$, i.e., $i \equiv m\equiv 0~(mod~m)$. Thus $i=0$. This implies that $\mathsf{Stab}_H (A_0) = \{\mathsf{id}\}$ and hence the theorem holds.\qed

{\theorem \label{family-3-non-cayley-theorem} $R_{2m}(m-2,m-1)$ is not a Cayley graph, if $m$ is odd and $m\not\equiv 0~(mod~3)$.}\\
\pf Consider $K=\langle \varepsilon_0,\varepsilon_1,\ldots,\varepsilon_{m-1} \rangle$. Then $K\cong \mathbb{Z}_2^m$ and $|K|=2^m$ as $\circ(\varepsilon_i)=2$ and $\varepsilon_i\varepsilon_j=\varepsilon_j\varepsilon_i, \forall i,j \in \{0,1,\ldots,m-1\}$.

If possible, let $H$ be a regular subgroup of $G$. Then $|H|=4m$. Thus $$|HK|=\dfrac{|H||K|}{|H \cap K|}=\dfrac{2^2m\cdot 2^m}{|H \cap K|}\leq 2^{m+1}m, \mbox{ i.e., }|H \cap K|\geq 2.$$ 
Now,  as $|H|=4m$,  where $m$ is odd and $|K|=2^m$, we have $|H\cap K|=2$ or $4$. We will prove that $|H\cap K|=4$. In fact, using the next two claims, we prove that $|H\cap K|\neq 2$. \\
{\it Claim 1:} If $|H\cap K|=2$, then the non-identity element of $H\cap K$ must be $\varepsilon_0\varepsilon_1\cdots\varepsilon_{m-1}=\rho^m$.\\
{\it Proof of Claim 1:} Let $\alpha=\varepsilon_{l_1}\varepsilon_{l_2}\cdots\varepsilon_{l_p}$ be the non-identity element of $H\cap K$. Let $L=\langle \mu, \varepsilon_0,\varepsilon_1,\ldots,\varepsilon_{m-1} \rangle$. Then $|L|=2^{m+1}$ and $K \subsetneq L$ as $\mu \in L\setminus K$. Thus $$|HL|=\dfrac{|H||L|}{|H \cap L|}=\dfrac{4m\cdot 2^{m+1}}{|H \cap L|}\leq |G|=2^{m+1}m, \mbox{ i.e., }|H \cap L|\geq 4.$$
As $|H\cap K|=2$ and $K \subsetneq L$,  there exists atleast one element of the form $\beta=\mu \varepsilon_{i_1}\varepsilon_{i_2}\cdots\varepsilon_{i_s}$ in $H\cap L$.

Again, let $L'=\langle \rho\mu, \varepsilon_0,\varepsilon_1,\ldots,\varepsilon_{m-1} \rangle$. By similar arguments, we can deduce that $|H\cap L'|\geq 4$. So there exists an element of the form $\gamma=\rho\mu \varepsilon_{j_1}\varepsilon_{j_2}\cdots\varepsilon_{j_t}$ in $H\cap L'$.

As $\alpha,\beta,\gamma \in H$, it follows that $\beta\alpha\beta^{-1},\gamma\alpha\gamma^{-1}\in H$. Observe that $$\beta\alpha\beta^{-1}=(\mu \varepsilon_{i_1}\varepsilon_{i_2}\cdots\varepsilon_{i_s})(\varepsilon_{l_1}\varepsilon_{l_2}\cdots\varepsilon_{l_p})(\mu \varepsilon_{i_1}\varepsilon_{i_2}\cdots\varepsilon_{i_s})^{-1}=\mu(\varepsilon_{l_1}\varepsilon_{l_2}\cdots\varepsilon_{l_p})\mu.$$
As $\mu\varepsilon_i=\varepsilon_{m-1-i}\mu$, $\beta\alpha\beta^{-1}$ is product of some $\varepsilon_i$'s and hence $\mathsf{id}\neq \beta\alpha\beta^{-1}\in H\cap K$. Since $|H \cap  K|=2$, then $\alpha=\beta\alpha\beta^{-1}$.

Similarly, $$\gamma\alpha\gamma^{-1}=(\rho\mu \varepsilon_{j_1}\varepsilon_{j_2}\cdots\varepsilon_{j_t})(\varepsilon_{l_1}\varepsilon_{l_2}\cdots\varepsilon_{l_p})(\rho\mu \varepsilon_{j_1}\varepsilon_{j_2}\cdots\varepsilon_{j_t})^{-1}=\rho(\mu \varepsilon_{l_1}\varepsilon_{l_2}\cdots\varepsilon_{l_p}\mu)\rho^{-1}$$
$$=\rho(\beta\alpha\beta^{-1})\rho^{-1}=\rho\alpha\rho^{-1}.$$
As $\rho\varepsilon_i=\varepsilon_{i+1}\rho$, $\rho\alpha\rho^{-1}$ is product of some $\varepsilon_i$'s and hence $\gamma\alpha\gamma^{-1} \in H\cap K$ and by similar arguments, we have $\alpha=\gamma\alpha\gamma^{-1}$.

Thus, using $\rho\varepsilon_i=\varepsilon_{i+1}\rho$, we get
\begin{equation}\label{family-3-non-cayley-theorem-equation} \varepsilon_{l_1}\varepsilon_{l_2}\cdots\varepsilon_{l_p}=\alpha=\rho\alpha\rho^{-1}=\rho(\varepsilon_{l_1}\varepsilon_{l_2}\cdots\varepsilon_{l_p})\rho^{-1}=\varepsilon_{l_1+1}\varepsilon_{l_2+1}\cdots\varepsilon_{l_p+1}
\end{equation}

As $K=\langle \varepsilon_0,\varepsilon_1,\ldots,\varepsilon_{m-1}\rangle \cong \mathbb{Z}_2^m$ and $\varepsilon_i$'s corresponds to the standard generators of $\mathbb{Z}_2^m$, i.e., $\varepsilon_i\leftrightarrow (0,0,\ldots,0,1,0,\ldots,0)$ with the only $1$ occuring in the $(i+1)$th position, $\varepsilon_{l_1}\varepsilon_{l_2}\cdots\varepsilon_{l_p}$ corresponds to the vector in $\mathbb{Z}_2^m$ with $1$'s in $l_1+1,l_2+1,\ldots,l_p+1$ positions and $\varepsilon_{l_1+1}\varepsilon_{l_2+1}\cdots\varepsilon_{l_p+1}$ corresponds to the vector with $1$'s in $l_1+2,l_2+2,\ldots,l_p+2$ positions. Thus, from Equation \ref{family-3-non-cayley-theorem-equation}, we get that all the positions in the vector must be $1$, i.e., $\alpha=\varepsilon_0\varepsilon_1\cdots\varepsilon_{m-1}=\rho^m$. Hence the claim is true.

\noindent {\it Claim 2:} If $|H\cap K|=2$, then $\rho^m\not\in H$\\
{\it Proof of Claim 2:} As $H \cap L$ is a subgroup of $H$ and $m$ is odd, therefore $4\leq |H \cap L|\mid 4m$ implies $|H \cap L|=4$. Thus $H\cap L$ is either isomorphic to $\mathbb{Z}_2\times \mathbb{Z}_2$ or $\mathbb{Z}_4$. Note that any non-identity element $\sigma \in H\cap L$ must contain in its expression either $\varepsilon_0$ or $\varepsilon_{m-1}$, as otherwise $\sigma \in \langle \mu,\varepsilon_1,\varepsilon_2,\ldots,\varepsilon_{m-2} \rangle=\mathsf{Stab}_G(A_0)$, a contradiction to the fact that $\sigma$ belongs to a regular subgroup $H$.

 Suppose that $H\cap L$ is isomorphic to $\mathbb{Z}_2\times \mathbb{Z}_2$. As $ H\cap K \subsetneq H \cap L$, therefore there exists a non-identity element in $H\cap L$ of the form $\sigma=\mu\varepsilon_{i_1}\varepsilon_{i_2}\cdots \varepsilon_{i_s}$. As explained earlier, $\sigma$ must contain in its expression either $\varepsilon_0$ or $\varepsilon_{m-1}$. In fact, in this case, both $\varepsilon_0$ and $\varepsilon_{m-1}$ must occur in the expression of $\sigma$, as otherwise $\circ(\sigma)=4$. Note that by Claim 1, $\rho^m \in H\cap L$. Thus, for all the three non-identity elements, $\rho^m, \sigma,\sigma'$ (say) in $H\cap L$, both $\varepsilon_0$ and $\varepsilon_{m-1}$ must occur. Also as $H\cap L\cong \mathbb{Z}_2\times \mathbb{Z}_2$, we have $\sigma\sigma'=\rho^m$. But if $\sigma,\sigma'$ contains both $\varepsilon_0$ and $\varepsilon_{m-1}$, then $\rho^m$ contains neither $\varepsilon_0$ nor $\varepsilon_{m-1}$, a contradiction. Hence $H\cap L\not\cong \mathbb{Z}_2\times \mathbb{Z}_2$.
 
 Suppose that $H\cap L$ is isomorphic to $\mathbb{Z}_4$. As $\circ(\rho^m)=2$, there exists a non-identity element $\zeta=\mu\varepsilon_{j_1}\varepsilon_{j_2}\cdots \varepsilon_{j_s} \in H\cap L$ such that $\langle \zeta \rangle=H\cap L$ and $\zeta^2=\rho^m$. Note that the number of $\varepsilon_i$'s in the expression of $\zeta^2$ is always even but that of $\rho^m$ is $m$ (odd) as $\rho^m=\varepsilon_0\varepsilon_1\cdots\varepsilon_{m-1}$. Hence, $H\cap L\not\cong \mathbb{Z}_4$.
 
 Thus, by Claim 1 and 2, we get $|H\cap K|=4$.  As $K \cong \mathbb{Z}_2^m$, we have $H\cap K \cong \mathbb{Z}_2 \times \mathbb{Z}_2$. Recall that $$ \mathsf{Stab}_G(B_{(m+3)/2})=\langle \rho\mu, \varepsilon_0,\varepsilon_1,\ldots,\varepsilon_{(m+1)/2},\varepsilon_{(m+7)/2},\ldots, \varepsilon_{m-1} \rangle.$$ Again, as the graph is vertex-transitive, by orbit-stabilizer theorem, we have $G=H\cdot \mathsf{Stab}_G(B_{(m+3)/2})$.  Thus, $\rho=hb$, where $h\in H$ and $b \in \mathsf{Stab}_G(B_{(m+3)/2})$. 
 
 \noindent {\it Claim 3:} $\rho\mu$ does not occur in the expression of $b$.\\
 {\it Proof of Claim 3:} If possible, let $b=\rho\mu\varepsilon_{l_1}\varepsilon_{l_2}\cdots\varepsilon_{l_p}$ and hence $h=\rho b^{-1}=\mu\varepsilon_{t_1}\varepsilon_{t_2}\cdots\varepsilon_{t_p}\in H \cap L$. Again, as $H\cap K \subseteq H \cap L$ and $|H\cap L|=|H\cap K|=4$, we have $H\cap K=H\cap L$. Thus, $h \in H\cap K \subset K$ and hence $h$ does not contain $\mu$ in its expression, a contradiction. Thus Claim 3 is true.

Therefore, by Claim 3, $b=\varepsilon_{l_1}\varepsilon_{l_2}\cdots\varepsilon_{l_p}$ and $h=\rho b^{-1}=\rho \varepsilon_{l_1}\varepsilon_{l_2}\cdots\varepsilon_{l_p} \in H$.

Let $H \cap K=\{\mathsf{id},\alpha_1,\alpha_2,\alpha_3\}\cong \mathbb{Z}_2 \times \mathbb{Z}_2$. Thus $h\alpha_i h^{-1}\in H$.  As $\alpha_i$'s, being elements of $K$, are product of some $ \varepsilon_i$'s and $\varepsilon_i \varepsilon_j=\varepsilon_j \varepsilon_i$, $\rho\varepsilon_i=\varepsilon_{i+1}\rho$, we have  
\begin{equation}\label{family-3-non-cayley-claim-4-equation}
h\alpha_i h^{-1}=\rho \alpha_i \rho^{-1}=\rho( \varepsilon_{i_1}\varepsilon_{i_2}\cdots\varepsilon_{i_s}) \rho^{-1}=\varepsilon_{i_1+1}\varepsilon_{i_2+1}\cdots\varepsilon_{i_s+1}\in K \mbox{ for }i=1,2,3. 
\end{equation}
Thus $h\alpha_i h^{-1} \in H \cap K=\{\mathsf{id},\alpha_1,\alpha_2,\alpha_3\}$.

\noindent{\it Claim 4:} $h\alpha_1 h^{-1}=\alpha_2$ or $\alpha_3$.\\
{\it Proof  of Claim 4:} If $h\alpha_1 h^{-1}=\mathsf{id}$, then $\alpha_1=\mathsf{id}$, a contradiction.

 If $h\alpha_1 h^{-1}=\alpha_1$, then as above, get $\varepsilon_{i_1+1}\varepsilon_{i_2+1}\cdots\varepsilon_{i_s+1}=\varepsilon_{i_1}\varepsilon_{i_2}\cdots\varepsilon_{i_s}.$ Now, as in proof of Claim 1, we can argue that this implies $\alpha_1=\rho^m$. But, in that case, we must have $h\alpha_2 h^{-1}=\alpha_3$ and $h\alpha_3 h^{-1}=\alpha_2$,  because otherwise 
\begin{itemize}
	\item $h\alpha_2 h^{-1}=\mathsf{id}$ implies $\alpha_1=\mathsf{id}$, a contradiction.
	\item $h\alpha_2 h^{-1}=\mathsf{\alpha_2}$ implies $\alpha_2=\rho^m$, a contradiction, as $\alpha_1\neq \alpha_2$.
	\item $h\alpha_2 h^{-1}=\mathsf{\alpha_1}$ implies $h\alpha_2 h^{-1}=h\alpha_1 h^{-1}$, i.e., $\alpha_1=\alpha_2$, a contradiction.
\end{itemize}
Thus we have $h\alpha_2 h^{-1}=\rho\alpha_2 \rho^{-1}=\alpha_3$ and $h\alpha_3 h^{-1}=\rho\alpha_3 \rho^{-1}=\alpha_2$. Hence, from Equation \ref{family-3-non-cayley-claim-4-equation}, we see that both $\alpha_2$ and $\alpha_3$ are product of $\varepsilon_i$'s and the number of $\varepsilon_i$'s occuring in their expressions are same. Thus the number of $\varepsilon_i$'s occuring in the expression of $\alpha_2\alpha_3$ is even. However, $\alpha_2\alpha_3=\alpha=\rho^m=\varepsilon_0\varepsilon_1\cdots\varepsilon_{m-1}$ has odd number of $\varepsilon_i$'s occuring in its expression. This is a contradiction and hence $h\alpha_1 h^{-1}\neq \alpha_1$. Thus Claim 4 is true.

Without loss of generality, we can assume that $h\alpha_1 h^{-1}=\alpha_2$. Thus $h\alpha_2 h^{-1}$ is either $\alpha_1$ or $\alpha_3$. If $h\alpha_2 h^{-1}=\alpha_1$, we must have $h\alpha_3 h^{-1}=\alpha_3$, a contradiction, as shown in Claim 4. Hence we have $h\alpha_2 h^{-1}=\alpha_3$ and similarly $h\alpha_3 h^{-1}=\alpha_1$. So, by Equation \ref{family-3-non-cayley-claim-4-equation}, we get $\rho\alpha_1\rho^{-1}=\alpha_2$, $\rho\alpha_2\rho^{-1}=\alpha_3$ and $\rho\alpha_3\rho^{-1}=\alpha_1$. Hence, we have $$\alpha_1=\rho\alpha_3\rho^{-1}=\rho(\rho\alpha_2\rho^{-1})\rho^{-1}=\rho^2(\rho\alpha_1\rho^{-1})\rho^{-2}=\rho^3\alpha_1\rho^{-3}, ~i.e.,~ \rho^3\alpha_1=\alpha_1\rho^3.$$
Similarly, we have $\rho^3\alpha_2=\alpha_2\rho^3$ and $\rho^3\alpha_3=\alpha_3\rho^3$.

Recall that $H \cap K=\{\mathsf{id},\alpha_1,\alpha_2,\alpha_3\}\cong \mathbb{Z}_2 \times \mathbb{Z}_2$ and $\alpha_i$'s are product of some $ \varepsilon_j$'s. Let $$\alpha_1=\varepsilon_{i_1}\varepsilon_{i_2}\cdots \varepsilon_{i_l}; \alpha_2=\varepsilon_{j_1}\varepsilon_{j_2}\cdots \varepsilon_{j_p}; \alpha_3=\varepsilon_{k_1}\varepsilon_{k_2}\cdots \varepsilon_{k_q}.$$ Note that each $\alpha_i$ must contain either $\varepsilon_0$ or $\varepsilon_{m-1}$ in its expression, as otherwise it will be an element of $\mathsf{Stab}_G(A_0)$ and hence can not belong to $H$. As $\alpha_1\alpha_2=\alpha_3$ and $\alpha_1\alpha_2\alpha_3=\mathsf{id}$, without loss of generality, we can assume that, among $\varepsilon_0$ or $\varepsilon_{m-1}$, $\alpha_1$ contains only $\varepsilon_0$, $\alpha_2$ contains only $\varepsilon_{m-1}$ and $\alpha_3$ contains both $\varepsilon_0$ and $\varepsilon_{m-1}$ in their expressions. This happens because if two of the $\alpha_i$'s contain both $\varepsilon_0$ and $\varepsilon_{m-1}$ in their expressions, then the their product, i.e., the third $\alpha_i$, will not have $\varepsilon_0$ or $\varepsilon_{m-1}$ in its expression, thereby making it an element of $\mathsf{Stab}_G(A_0)$.

Now, from the relation $\rho^3\alpha_1=\alpha_1\rho^3$ and using the fact that $\rho\varepsilon_i=\varepsilon_{i+1}\rho$, we get, $$(\varepsilon_{i_1}\varepsilon_{i_2}\cdots \varepsilon_{i_l})\rho^3=\rho^3(\varepsilon_{i_1}\varepsilon_{i_2}\cdots \varepsilon_{i_l})=(\varepsilon_{i_1+3}\varepsilon_{i_2+3}\cdots \varepsilon_{i_l+3})\rho^3,$$ 
$$ \mbox{ i.e., }\varepsilon_{i_1}\varepsilon_{i_2}\cdots \varepsilon_{i_l}=\varepsilon_{i_1+3}\varepsilon_{i_2+3}\cdots \varepsilon_{i_l+3}.$$

Now, as $m$ is not a multiple of $3$, $m$ is of the form $3t+1$ or $3t+2$. 

If $m=3t+1$, then by using the standard generators of $\mathbb{Z}_2^m$, as in the proof of Claim 1, we get that all of $\varepsilon_0,\varepsilon_3,\varepsilon_6,\ldots,\varepsilon_{3t}=\varepsilon_{m-1}$ occurs in the expression of $\alpha_1$, a contradiction to that fact that among $\varepsilon_0$ or $\varepsilon_{m-1}$, $\alpha_1$ contains only $\varepsilon_0$.

Similarly, if $m=3t+2$, we get all of $$\varepsilon_0,\varepsilon_3,\varepsilon_6,\ldots,\varepsilon_{3t}=\varepsilon_{m-2},\varepsilon_1,\varepsilon_4,\cdots,\varepsilon_{3t+1}=\varepsilon_{m-1}$$ occurs in the expression of $\alpha_1$, a contradiction.

Thus, we conclude that there does not exist any regular subgroup $H$ of $\mathsf{Aut}(R_{2m}(m-2,m-1))$ and hence $R_{2m}(m-2,m-1)$ is not a Cayley graph, when $m$ is odd and not a multiple of $3$. 
  \qed

\section{Family-4 [$R_{12m}(\pm (3m+2),\pm (3m-1))$ and $R_{12m}(\pm (3m-2),\pm (3m+1))$]}
As $R_n(a,r)=R_n(a,-r)$ and $R_n(a,r)\cong R_n(-a,r)$, it is enough to check $R_{12m}( 3m+2,3m-1)$ and $R_{12m}(3m-2,3m+1)$. More precisely, it suffices to work with the family $R_{12m}(3d+2,9d+1)$ where $d=\pm m~(mod~12m)$, as mentioned in Section 3.3 of \cite{rotary}.
Define $\sigma$ as follows:
$$\sigma(A_i)=\left\lbrace \begin{array}{lr}
A_i & \mbox{ if }i\equiv 0~(mod~3)\\
B_{i-1} & \mbox{ if }i\equiv 1~(mod~3)\\
B_{i-1-3d} & \mbox{ if }i\equiv 2~(mod~3)
\end{array} \right. \mbox{ and }~\sigma(B_i)=\left\lbrace \begin{array}{lr}
A_{i+1} & \mbox{ if }i\equiv 0~(mod~3)\\
A_{i+3d+1} & \mbox{ if }i\equiv 1~(mod~3)\\
B_{i+6d} & \mbox{ if }i\equiv 2~(mod~3)
\end{array} \right.$$
Also, if $m\equiv 2~(mod~4)$, let $b=d+1$ and define $\omega$ as follows:
$$\omega(A_i)=\left\lbrace \begin{array}{lr}
A_{bi} & \mbox{ if }i\equiv 0~(mod~3)\\
B_{bi-b} & \mbox{ if }i\equiv 1~(mod~3)\\
B_{b+bi-1} & \mbox{ if }i\equiv 2~(mod~3)
\end{array} \right. \mbox{ and }~\omega(B_i)=\left\lbrace \begin{array}{lr}
A_{bi+1} & \mbox{ if }i\equiv 0~(mod~3)\\
A_{4+bi-4b} & \mbox{ if }i\equiv 1~(mod~3)\\
B_{b+bi-1} & \mbox{ if }i\equiv 2~(mod~3)
\end{array} \right.$$
It was shown in \cite{rotary}, that $$G:=\mathsf{Aut}(R_{12m}(3d+2,9d+1))=\left\lbrace 
\begin{array}{lr}
\langle \rho, \mu,\sigma,\omega \rangle, & \mbox{ if }m\equiv 2~(mod~4)\\
\langle \rho, \mu,\sigma \rangle, & \mbox{ otherwise}
\end{array}
\right.$$
It is to be noted that $m\equiv 2~(mod~4)$ if and only if $-m\equiv 2~(mod~4)$. Thus, it is enough to work only with the family $R_{12m}(3m+2,9m+1)$.

{\theorem \label{family-4-m-odd-cayley-theorem} If $m$ is odd and $m\neq 3$, then $R_{12m}(3m+2,9m+1)$ is a Cayley graph.}\\
\pf As $m$ is odd, $G=\langle \rho, \mu,\sigma \rangle$. It can also be checked that $\sigma\rho^3\sigma=\rho^3; \sigma\mu=\mu\sigma; (\rho\sigma)^3=\rho^{3(m+1)}; \circ (\sigma)=2$. Let $\alpha=(\rho\sigma)^2$ and $\beta=\rho^2\mu\sigma$. As $m$ is odd and $m\neq 3$, it can be shown that $\circ (\alpha)=3m, \circ(\beta)=8$ and $\beta\alpha=\alpha^{-1}\beta^{-1}$. Define $$H=\langle \alpha,\beta: \circ (\alpha)=3m, \circ(\beta)=8; \beta\alpha=\alpha^{-1}\beta^{-1} \rangle$$
$$=\{\alpha^i\beta^j:0\leq i\leq 3m-1; 0\leq j \leq 7 \}~~~~~~$$
\noindent {\it Claim 1:} The elements in $H$ are distinct. \\ If not, suppose $$\alpha^{i_1}\beta^{j_1}=\alpha^{i_2}\beta^{j_2}, \mbox{ where } 0\leq i_1,i_2 < 3m, 0\leq j_1,j_2\leq 8,$$ i.e., 
\begin{equation}\label{alpha-beta-equation}
\alpha^{i_1-i_2}=\beta^{j_2-j_1}.
\end{equation}
$$\mbox{As } \alpha(A_0)=B_1,\alpha^2(A_0)=A_{3m+4},\alpha^3(A_0)=A_{6m+6},\alpha^4(A_0)=A_{6m+7},\ldots,\alpha^{3m}(A_0)=A_0,$$ any power of $\alpha$ maps $A_0$ to $A_{0(mod~3)}$ or $A_{1(mod~3)}$ or $B_{1(mod~3)}$. On the other hand, as $$\beta(A_0)=A_2,\beta^2(A_0)=B_{3m-1},\beta^3(A_0)=B_{3m+1},\beta^4(A_0)=A_{6m},$$
$$\beta^5(A_0)=A_{6m+2},\beta^6(A_0)=B_{6m-1},\beta^7(A_0)=B_{9m+1},\beta^8(A_0)=A_0,$$ we see that $\beta,\beta^2,\beta^5$ and $\beta^6$ maps $A_0$ to $A_{2(mod~3)}$. Thus, $j_2-j_1$ in Equation \ref{alpha-beta-equation} can take values from $\{0,3,4,7\}$.\\
If $j_2-j_1=0$, then it is obvious that $i_1=i_2$ and $j_1=j_2$.\\
If $j_2-j_1=4$, squaring Equation \ref{alpha-beta-equation}, we get, $\alpha^{2(i_1-i_2)}=\mathsf{id}$. Therefore, $3m|2(i_1-i_2)$. Now, as $gcd(2,3)=1$ and $m$ is odd, we have $3m|(i_1-i_2)$, i.e., $i_1=i_2$ and hence $j_1=j_2$.\\
If $j_2-j_1=3$, since $gcd(3,8)=1$, then $\circ (\beta^{j_2-j_1})=8$. Therefore, $\alpha^{8(i_1-i_2)}=\mathsf{id}$, i.e., $3m|8(i_1-i_2)$. As $m$ is odd, $3m$ is coprime to $8$ and hence, $3m|(i_1-i_2)$, i.e., $i_1=i_2$ and $j_1=j_2$.\\
The case $j_2-j_1=7$ follows similarly as above. Thus combining all the cases, we see that elements of $H$ are distinct and $H=3m\times 8=24m$.

\noindent {\it Claim 2:} $H$ acts transitively on $R_{12m}(3m+2,9m+1)$. \\In order to prove it, we show that the orbit of $A_0$, $\mathcal{O}_{A_0}$, under the action of $H$ is the vertex set of $R_{12m}(3m+2,9m+1)$. By orbit-stabilizer theorem, we get $$|\mathcal{O}_{A_0}|=\dfrac{|H|}{|\mathsf{Stab}_H(A_0)|}.$$ As the number of vertices in $R_{12m}(3m+2,9m+1)$ is $24m$ and $|H|=24m$, it is enough to show that $\mathsf{Stab}_H(A_0)=\{\mathsf{id}\}$. Let $\alpha^i\beta^j$ be an arbitrary element of $H$ which stabilizes $A_0$, i.e., $\alpha^{-i}(A_0)=\beta^j(A_0)$ with $0\leq i\leq 3m-1; 0\leq j \leq 7$.  Again, by mimicing the argument used in the proof of {\it Claim 1}, one can conclude that $j\in \{0,3,4,7\}$.\\
If $j=4$, then $\alpha^{-i}(A_0)=\beta^4(A_0)=A_{6m}$. Thus, $-i$ and hence $i$ is a multiple of $3$. [since, $\alpha^x$ sends $A_0$ to $A_{0(mod~3)}$, only if $x$ is a multiple of $3$] Let $-i=3k$ and therefore $A_{6m}=\alpha^{3k}(A_0)=A_{k(6m+6)}$, i.e., $12m|k(6m+6)-6m$, i.e., $2m|m(k-1)+k$, i.e., $m|k$ which implies $k=lm$. Again, as $2m|m(k-1)+lm$, we have $2|k-1+l$, i.e., $2|l(m+1)-1$. But this is a contradiction, as $m+1$ is even and hence $l(m+1)-1$ is odd. Thus $j\neq 4$.\\
If $j=3$, then $\alpha^{-i}(A_0)=\beta^3(A_0)=B_{3m+1}$.  As $3m+1\equiv 1(mod~3)$, we have $-i=3k+1$ [since, $\alpha^x$ sends $A_0$ to $B_{1(mod~3)}$, only if $x\equiv 1(mod~3)$] Therefore, $\beta^3(A_0)=B_{3m+1}=\alpha^{3k+1}(A_0)=\alpha^{3k}(B_1)$, i.e., $B_{3m+1}=B_{1+6mk+6k}$. This implies $12m|6mk+6k-3m$, i.e., $4m|2mk+2k-m$,i.e., $m|2k$ and, as $m$ is odd, we have $m|k$. Let $k=lm$. Again, as $4m|2mk+2lm-m$, we have $4|2k+2l-1$. However, this is a contradiction, as $2k+2l-1$ is odd and hence $j \neq 3$.\\
Using similar arguments as above, it can be shown that $j\neq 7$.\\
Thus, we have $j=0$ and this, in turn, implies $i=0$. Hence, $\mathsf{Stab}_H(A_0)=\{\mathsf{id}\}$.

Finally, in view of Remark \ref{regular-transitive}, $H$ acts regularly on $R_{12m}(3m+2,9m+1)$ and hence $R_{12m}(3m+2,9m+1)$ is a Cayley graph, if $m$ is odd and $m\neq 3$.\qed

{\theorem \label{family-4-m=3-theorem} If $m= 3$, then $R_{12m}(3m+2,9m+1)$, i.e, $R_{36}(11,28)$ is a Cayley graph.}\\
\pf This can be checked using Sage programming. See Appendix for the SageMath code.\qed

{\theorem If $m\equiv 0(mod~4)$, then $R_{12m}(3m+2,9m+1)$ is not a Cayley graph.}\\
\pf As $m\not\equiv 2(mod~4)$, $$G=\langle \rho,\mu,\sigma: \rho^n=\mu^2=\sigma^2=\mathsf{id}; \mu\rho\mu=\rho^{-1}, \sigma\rho^3\sigma=\rho^3,\sigma\mu=\mu\sigma,~~~~~~~~~~~~~~~~~~$$ 
$$(\rho\sigma)^3=\rho^{3(m+1)}, (\rho\sigma\rho)^3=\rho^{9m+6} \rangle, \mbox{ where }n=12m$$
If possible, let $R_{12m}(3m+2,9m+1)$ be a Cayley graph, $H$ be a regular subgroup of $G$ and $K=\mathsf{Stab}_G(A_0)$. Then $|G|=96m=8n$ (See Lemma \ref{family-4-calc} in Appendix), $|H|=2n=24m$ and $H\cap K=\{\mathsf{id}\}$. 

Let $K'=\langle \rho \rangle$. Then $|K'|=n$ and $|HK'|=\frac{|H||K'|}{|H\cap K'|}=\frac{2n^2}{n/t}\leq |G|=8n$, where $t$ is a factor of $n$. Thus, $t\leq 4$, i.e., $t=1,2,3$ or $4$. If $t=1$, then $H\cap K'=K'$, i.e., $\rho \in H$. If $t=2$, then $H\cap K'=\langle \rho^2 \rangle$, i.e., $\rho^2 \in H$. If $t=3$, then $H\cap K'=\langle \rho^3 \rangle$, i.e., $\rho^3 \in H$. If $t=4$, then $H\cap K'=\langle \rho^4 \rangle$, i.e., $\rho^4 \in H$. Combining all the cases, we get that 
\begin{equation}\label{family-4-not-cayley-equation}
\mbox{either }\rho^3 \in H \mbox{ or }\rho^4 \in H.
\end{equation}
{\it Claim:}  $\rho^4 \in H$.\\
{\it Proof of Claim:} Suppose that that $\rho^3\in H$ but $\rho^4 \not\in H$. Let $L=\langle \rho,\mu \rangle$. Then $|L|=2n$. Therefore $$|HL|=\dfrac{|H||L|}{|H \cap L|}=\dfrac{2n\cdot 2n}{2n/t}=2nt\leq |G|=8n, \mbox{ i.e., }t=1,2,3 \mbox{ or }4 \mbox{ and } t \mbox{ divides }2n.$$
Therefore, $|H\cap L|=2n,n,2n/3$ or $n/2$, i.e., $|H\cap L|\geq n/2$. As $\rho^3 \in H \cap L$, we have $\langle \rho^3 \rangle\subseteq H\cap L$ and $|\langle \rho^3 \rangle|=n/3$. Thus, $(H\cap L)\setminus  \langle \rho^3 \rangle \neq \emptyset$. 

Now, as $\rho^{2i}\mu(A_i)=A_i$, $\rho^{2i}\mu \not\in H$. Similarly, if $\rho^{2i+1}\mu \in H$, then $\rho^3\cdot \rho^{2i+1}\mu \in H$, i.e., $\rho^{2i+4}\mu \in H$. Note that $2i+4$ is even and hence by previous argument, $\rho^{2i+4}\mu \not\in H$, i.e., $\rho^{2i+1}\mu \not\in H$. This shows that $H$ does not contain any element of the form $\rho^i\mu$. Moreover, $\mu \not\in H$. Now, as $(H\cap L)\setminus  \langle \rho^3 \rangle \neq \emptyset$, $H$ must contain an element of the form $\rho^i$, where $i$ is not a multiple of $3$. Again, as $\rho^3 \in H$, either $\rho$ or $\rho^2 \in H$, i.e., $\rho^4 \in H$. This is a contradiction to the assumption that $\rho^4 \not\in H$. Thus the claim is true.

Let $K''=\langle \rho\sigma\rangle$. As $\circ(\rho\sigma)=n$, we have $|K''|=n$ and by similar arguments as above, we get that either $(\rho\sigma)^3\in H$ or $(\rho\sigma)^4 \in H$.

Case 1: If $\rho^4 \in H$ and $(\rho\sigma)^4 \in H$, then $$(\rho\sigma)^4=(\rho\sigma)^3(\rho\sigma)=\rho^{3(m+1)}\rho\sigma=\rho^{3m+4}\sigma=\rho^{12l+4}\sigma=(\rho^4)^{3l+1}\sigma\in H ~[\mbox{letting }m=4l].$$
As $\rho^4 \in H$, therefore $\sigma \in H$. But as $\sigma(A_0)=A_0$, i.e., $\sigma$ stabilizes $A_0$, it can not be in $H$. This is a contradiction.

Case 2: If $\rho^4 \in H$ and $(\rho\sigma)^3 \in H$, then $(\rho\sigma)^3=\rho^{3(m+1)}=\rho^{12l+3}=(\rho^4)^{3l}\rho^3 \in H$, where $m=4l$ i.e., $\rho^3 \in H$. Again, as $\rho^4 \in H$, we have $\rho \in H$. As $\circ(\rho)=n$ and $[H:\langle \rho \rangle]=2$, $\langle \rho \rangle$ is normal in $H$.

From definition, it follows that $\mathsf{id}, \mu,\sigma,\mu\sigma \in K$. On the other hand, as $R_{12m}(3m+2,9m+1)$ is vertex transitive, by orbit-stabilizer theorem, we have $$|K|=\dfrac{|G|}{2n}=\dfrac{8n}{2n}=4. \mbox{ Hence, }K=\mathsf{Stab}_G(A_0)=\{\mathsf{id}, \mu,\sigma,\mu\sigma\} \mbox{ and }|HK|=\dfrac{2n\cdot 4}{1}=8n=|G|.$$
Thus, $HK=G$. As $\sigma\rho \in G$, it can be expressed in the form $\alpha\beta$, where $\alpha \in H$ and $\beta \in K=\{\mathsf{id}, \mu,\sigma,\mu\sigma\}$. \\
If $\beta=\mathsf{id}$, then $\alpha=\sigma\rho \in H$, i.e., $\sigma \in H$ (as $\rho \in H$), which is a contradiction, as $H$, being a regular subgroup can not contain any non-identity element which stabilizes $A_0$. \\
If $\beta=\mu$, then $\sigma\rho=\alpha\mu$, i.e., $\alpha=\sigma\mu\rho^{-1}\in H$, i.e., $\sigma\mu \in H$ (as $\rho \in H$), which is a contradiction. \\
If $\beta=\sigma$, then $\alpha=\sigma\rho\sigma \in H$. Since $\langle \rho \rangle$ is normal in $H$, therefore $(\sigma\rho\sigma)\rho(\sigma\rho\sigma)^{-1}\in H$, i.e., $$(\sigma\rho\sigma)\rho(\sigma\rho\sigma)^{-1}=(\sigma\rho\sigma)\rho\sigma\rho^{-1}\sigma=(\sigma\rho)^3\rho^{-2}\sigma=\rho^{3m+1}\sigma \in H\Rightarrow \sigma \in H~(\mbox{as }\rho \in H),$$ a contradiction.\\ 
If $\beta=\mu\sigma$, then $\sigma\rho=\alpha\mu\sigma$, i.e., $\alpha=\sigma\rho\mu\sigma\in H$. Since $\langle \rho \rangle$ is normal in $H$, therefore $(\sigma\rho\mu\sigma)\rho(\sigma\rho\mu\sigma)^{-1}\in H$, i.e.,
$$(\sigma\rho\mu\sigma)\rho(\sigma\mu\rho^{-1}\sigma)=(\sigma\rho\mu\sigma)\rho(\sigma\rho\mu\sigma)=\sigma\rho\mu(\sigma\rho)^2\mu\sigma~~~~~~~~~~~~~~~~~~~~~~~~~~~~~~~~~~~~~~~~~~~$$
$$~~~~~~~~~~~~~~~~~~~~~~=\sigma\rho\mu(\rho^{3m+2}\sigma)\mu\sigma~~~~[\mbox{as }(\sigma\rho)^3=\rho^{3m+3}, \mbox{ we have }(\sigma\rho)^2=\rho^{3m+2}\sigma]$$
$$=\sigma\rho\mu\rho^{3m+2}\mu~~~~~~~~~~[\mbox{as }\sigma\mu=\mu\sigma \mbox{ and }\sigma^2=\mathsf{id}]~~$$
$$~~~~~~~~~~~~~~~~~~~~~~~~=\sigma\rho\rho^{-3m-2}=\sigma\rho^{-3m-1}\in H \Rightarrow \sigma \in H (\mbox{as } \rho \in H), \mbox{ a contradiction.}$$

Thus, combining Case 1 and 2, we conclude that there does not exist any regular subgroup $H$ of $G$, i.e., $R_{12m}(3m+2,9m+1)$ is not a Cayley graph, if $m\equiv 0(mod~4)$. \qed

\subsection{$m\equiv 2(mod~4)$}
 As $m\equiv 2(mod~4)$, $G=\langle \rho, \mu,\sigma,\omega \rangle$. It can be checked that $\sigma\rho^3\sigma=\rho^3; \sigma\mu=\mu\sigma; \sigma\omega=\omega\sigma; \omega\rho=\sigma\rho\omega; \omega\mu=\mu\omega\sigma; (\rho\sigma)^3=\rho^{3(m+1)}; \omega\rho^{3l}=\rho^{3l(m+1)}; (\rho\sigma\rho)^3=(\rho^2\sigma)^3=\rho^{9m+6};  \circ (\sigma)=\circ (\omega)=\circ (\sigma\omega)=2; \circ (\omega\mu)=4$.

Let $\alpha=\omega\sigma\rho^{4m}\omega\sigma$ and $\beta=\rho^{3m/2}$. Using the above relations,  it can be shown that $\circ(\alpha)=3; \circ(\beta)=8;\alpha\beta=\beta\alpha$. Define $$\gamma=\left\lbrace 
\begin{array}{ll}
\rho^{8m}\sigma\rho^2\omega, & \mbox{ if }m \mbox{ is of the form }12l+2 \mbox{ or }12l+6\\
(\rho^{8m}\sigma\rho^2\omega)^3, & \mbox{ if }m \mbox{ is of the form }12l+10.
\end{array}
\right.$$

In all the cases, it can be checked that $\circ(\gamma)=2m$, $\alpha\gamma=\gamma\alpha$ and $\gamma\beta=\beta^{m+1}\gamma$. It is to be noted that $\alpha=\omega\sigma\rho^{4m}\omega\sigma=(\omega\sigma\rho\omega\sigma)^{4m}=[\omega(\sigma\rho\omega)\sigma]^{4m}=(\omega(\omega\rho)\sigma)^{4m}=(\rho\sigma)^{4m}$.

{\proposition \label{family-4-gamma-results}
	\begin{enumerate}
		\item $\gamma^2=\left\lbrace 
		\begin{array}{ll}
		\rho^{4m+4}, & \mbox{ if }m \mbox{ is of the form }12l+2 \mbox{ or }12l+6\\
		\rho^{12}, & \mbox{ if }m \mbox{ is of the form }12l+10.
		\end{array}
		\right.$
		\item $\gamma^m=\left\lbrace 
		\begin{array}{ll}
			\alpha^2\beta^4, & \mbox{ if }m \mbox{ is of the form }12l+6\\
			\beta^4, & \mbox{ if }m \mbox{ is of the form }12l+2 \mbox{ or }12l+10.
		\end{array}
		\right.$
		
	\end{enumerate}
	
	 }
\pf See Appendix.

{\theorem \label{family-4-m=12l+2-Cayley-theorem} If $m\equiv 2(mod~12)$, then $R_{12m}(3m+2,9m+1)$ is a Cayley graph.}\\
\pf Let $m=12l+2$. Therefore $8m=96l+16$, i.e., $8m-4=12(8l+1)$. Then $\gamma^2=\rho^{4m+4}$.(by Proposotion \ref{family-4-gamma-results}) Define $$H=\langle \alpha,\beta,\gamma: \alpha^3=\beta^8=\gamma^{2m}=\mathsf{id};\alpha\beta=\beta\alpha, \alpha\gamma=\gamma\alpha,  \gamma\beta=\beta^{m+1}\gamma, \gamma^m=\beta^4 \rangle.$$
Thus, it is clear that every element of $H$ is of the form $\alpha^i\beta^j\gamma^k$ where $i=0,1,2$; $j=0,1,\ldots,7$ and $k=0,1,\ldots,m-1$. 

\noindent {\it Claim 1:} $H=\{\alpha^i\beta^j\gamma^k:i=0,1,2; j=0,1,\ldots,7; k=0,1,\ldots,m-1\}$.\\
{\it Proof of Claim 1:} If possible, let there exist $i_1,i_2\in \{0,1,2 \}, j_1,j_2 \in \{0,1,\ldots,7\}$ and $k_1,k_2 \in \{0,1,\ldots, m-1\}$, such that $\alpha^{i_1}\beta^{j_1}\gamma^{k_1}=\alpha^{i_2}\beta^{j_2}\gamma^{k_2}$. As $\alpha\beta=\beta\alpha$ and $\alpha\gamma=\gamma\alpha$, we have $$\alpha^{i_2-i_1}=\beta^{j_1-j_2}\gamma^{k_1-k_2}.$$

\noindent {\it Case 1:} $k_1-k_2$ is even.\\
As $\gamma^2=\rho^{4m+4}$ and $\beta=\rho^{3m/2}$, we have $\alpha^{i_2-i_1}=\rho^x$, i.e., $(\rho\sigma)^{4m(i_2-i_1)}=\rho^x$. This implies that $3$ divides $4m(i_2-i_1)$, i.e., $3|m$ or $3|(i_2-i_1)$. As $3$ does not divide $m$, we have $3|(i_2-i_1)$, i.e., $i_1=i_2$. Thus $\beta^{j_1-j_2}=\gamma^{k_2-k_1}=(\gamma^2)^{(k_2-k_1)/2}$, i.e., 
\begin{equation}\label{12k+2-equation-0}
(\rho)^{3m(j_1-j_2)/2}=(\rho^{4m+4})^{(k_2-k_1)/2}=\rho^{2(m+1)(k_2-k_1)}.
\end{equation}
Therefore, $12m$ divides $3m(j_1-j_2)/2 -2(m+1)(k_2-k_1)$, i.e.,
\begin{equation}\label{12k+2-equation-1}
24m \mbox{ divides }3m(j_1-j_2) -4(m+1)(k_2-k_1)
\end{equation}
Thus, $m$ divides $4(m+1)(k_2-k_1)$. As $gcd(m,4)=2$ and $gcd(d,d+1)=1$, it follows that $m/2$ divides $k_2-k_1$, i.e., $k_2-k_1=\frac{m}{2}s$. Since, $0\leq k_2-k_1<m$, we have $s=0$ or $1$. Again, as $m+1$ is a multiple of $3$, from Equation \ref{12k+2-equation-1}, we get that $12$ divides $3m(j_1-j_2)$, i.e., $2$ divides $j_1-j_2$. Let $j_1-j_2=2t$. As $0\leq j_1- j_2<8$, we have $t\in \{0,1,2,3\}$. Thus, rewriting Equation \ref{12k+2-equation-1}, we get $24m$ divides $6mt-2m(m+1)s$, i.e., $12$ divides $3t-(m+1)s$. Thus 
\begin{equation}\label{12k+2-equation-2}
4 \mbox{ divides }  \left( t-\frac{m+1}{3}s\right)=t-(4l+1)s, \mbox{ where }s\in \{0,1\}, t \in \{0,1,2,3\}.
\end{equation}
If $s=1$, then $k_2-k_1=m/2=6l+1$ is odd, a contradiction. Thus $s=0$ and hence from Equation \ref{12k+2-equation-2}, we have $4$ divides $t$, i.e., $t=0$. Therefore, we have $j_1-j_2=k_1-k_2=0$, and as a result $i_1=i_2$. Thus Claim 1 is true, if Case 1 holds. 

\noindent {\it Case 2:} $k_1-k_2$ is odd.\\
Let $k_1-k_2=2t-1$. Then we have $\alpha^{i_2-i_1}=\beta^{j_1-j_2}(\gamma^2)^{t}\gamma^{-1}.$ As $\gamma^2=\rho^{4m+4}$ and $\beta=\rho^{3m/2}$, we have $\gamma\alpha^{i_2-i_1}=\rho^x$. Now $i_2-i_1=0,1$ or $2$. Thus either of $\gamma,\alpha\gamma,\alpha^2\gamma$ is $\rho^x$. But $$\gamma(A_0)=\rho^{8m}\sigma\rho^2\omega(A_0)=\rho^{8m}\sigma(A_2)=\rho^{8m}(B_{9m+1})=B_{5m+1}~~~~~~~~~~~~~~~~~~~~~~~~~~~~~~~~~~~~~~~~~~$$
$$\alpha\gamma(B_0)=(\rho\sigma)^{4m}(A_{8m+3})=(\rho\sigma)^{48l+8}(A_{8m+3})=(\rho\sigma)^{2}((\rho\sigma)^{3})^{16l+2}(A_{8m+3})=A_{9m+3}~~~~~~~~~~~~~~~~$$
$$\alpha^2\gamma(A_0)=(\rho\sigma)^{8m}(B_{5m+1})=(\rho\sigma)((\rho\sigma)^{3})^{32l+5}(B_{5m+1})=(\rho\sigma)(B_{4m})=B_{10m+1}~~~~~~~~~~~~~~~~~~~~~~~~~~$$
As each of $\gamma,\alpha\gamma,\alpha^2\gamma$ maps some $A_i$ to some $B_j$, none of them is equal to $\rho^x$ and hence a contradiction. So $k_1-k_2$ can not be odd.\\ Combining Case 1 and 2, we conclude that Claim 1 is true and hence $|H|=24m=2n$. So, as in proof of Claim 2 in Theorem \ref{family-4-m-odd-cayley-theorem}, it suffices to show that $\mathsf{Stab}_H(A_0)=\{\mathsf{id}\}$. Let $\alpha^i\beta^j\gamma^k(A_0)=A_0$.\\
\noindent {\it Claim 2:} $k$ is even.\\
{\it Proof of Claim 2:} If possible, let $k$ be odd, say $k=2t+1$. Then, as $\alpha$ commutes with $\beta$ and $\gamma$, we have $\beta^j\gamma^{2t}\gamma\alpha^i(A_0)=A_0$, i.e., $\gamma\alpha^i(A_0)=\beta^{-j}(\gamma^2)^{-t}(A_0)=\rho^x(A_0)=A_x$, as in Case 2 above. Now, $i=0,1$ or $2$ and as $\gamma(A_0)=B_{5m+1}$ and $\alpha^2\gamma(A_0)=B_{10m+1}$, we have $i=1$. This implies $\alpha\beta^j\gamma^{2t+1}(A_0)=A_0$, i.e., $\beta^j(\gamma^2)^t\gamma(A_0)=\alpha^2(A_0)=A_{11m}$, i.e., $$A_{11m}=\beta^j(\gamma^2)^t\gamma(A_0)=\beta^j(\gamma^2)^t(B_{5m+1})=\rho^x(B_{5m+1})=B_{5m+x+1}, \mbox{ a contradiction.}$$ Hence the claim is true and let $k=2t$. Therefore, $$\beta^j(\gamma^2)^t(A_0)=\alpha^{-i}(A_0).$$ As left side of the above equation is $\rho^x(A_0)$ and $\alpha(A_0)=B_{10m-1}$, we conclude that $i=0$ or $1$. If $i=1$, then we have $\alpha\beta^j(\gamma^2)^t(A_0)=A_0$. Again as $\alpha$ commutes with $\beta$ and $\gamma$, we have $$A_0=\beta^j\gamma^{2t}\alpha(A_0)=\beta^j\gamma^{2t}(B_{10m-1})=\rho^x(B_{10m-1})=B_{10m+x-1}, \mbox{ a contradiction}.$$
Therefore, $i=0$ and hence we have $\beta^j(\gamma^2)^t(A_0)=A_0$, i.e., $$\rho^{4(m+1)t+3j\frac{m}{2}}(A_0)=A_0,\mbox{ i.e., }12m \mbox{ divides }4(m+1)t+3j\frac{m}{2}=12(4l+1)t+3j(6l+1)$$
Thus $12$ divides $3j(6l+1)$, i.e., $4|j(6l+1)$. However as $6l+1$ is odd and $j\in \{0,1,\ldots,7\}$, we have $j=0$ or $4$. If $j=4$, we have $12m$ divides $12(4l+1)t+12(6l+1)$, i.e., $m=12l+2=2(6l+1)$ divides $(4l+1)t+12(6l+1)$ and hence $2(6l+1)$ divides $(4l+1)t$. As $3(4l+1)-2(6l+1)=1$, we have $gcd(4l+1,6l+1)=1$ and hence $6l+1$ divides $t$. However as $0\leq k \leq m-1$, we have $0\leq t \leq \frac{m-1}{2}<6l+1$. Thus the only possible value of $t$ is $0$ and hence $k=0$. Therefore, we have $\beta^j(A_0)=A_0$, i.e., $\rho^{3(6l+1)j}(A_0)=A_0$. This implies that $12m=12(12l+2)$ divides $3(6l+1)j$, i.e., $8|j$ and hence $j=0$. 

Thus we have $\mathsf{Stab}_H(A_0)=\{\mathsf{id}\}$ and the theorem holds.\qed

{\theorem \label{family-4-m=12l+6-Cayley-theorem} If $m\equiv 6(mod~12)$, then $R_{12m}(3m+2,9m+1)$ is a Cayley graph.}\\
\pf Let $m=12l+6$. Therefore $8m=96l+48=12(8l+4)$. Also note that in this case, $\alpha=(\rho\sigma)^{4m}=((\rho\sigma)^3)^{4(4l+2)}=\rho^{12(m+1)(4l+2)}=\rho^{12(4l+2)}=\rho^{4m}$. Also $\gamma^2=\rho^{4m+4}$. (by Proposition \ref{family-4-gamma-results})
Define $$H=\langle \alpha,\beta,\gamma: \alpha^3=\beta^8=\gamma^{2m}=\mathsf{id};\alpha\beta=\beta\alpha, \alpha\gamma=\gamma\alpha,  \gamma\beta=\beta^{m+1}\gamma, \gamma^m=\alpha^2\beta^4 \rangle.$$
Thus, it is clear that every element of $H$ is of the form $\alpha^i\beta^j\gamma^k$ where $i=0,1,2$; $j=0,1,\ldots,7$ and $k=0,1,\ldots,m-1$. 

\noindent {\it Claim 1:} $H=\{\alpha^i\beta^j\gamma^k:i=0,1,2; j=0,1,\ldots,7; k=0,1,\ldots,m-1\}$.\\
{\it Proof of Claim 1:} If possible, let there exist $i_1,i_2\in \{0,1,2 \}, j_1,j_2 \in \{0,1,\ldots,7\}$ and $k_1,k_2 \in \{0,1,\ldots, m-1\}$, such that $\alpha^{i_1}\beta^{j_1}\gamma^{k_1}=\alpha^{i_2}\beta^{j_2}\gamma^{k_2}$. As $\alpha\beta=\beta\alpha$ and $\alpha\gamma=\gamma\alpha$, we have 
\begin{equation}\label{family-4-m=12l+6-claim-2-equation-1}
\alpha^{i_2-i_1}=\beta^{j_1-j_2}\gamma^{k_1-k_2}.
\end{equation}
If $k_1-k_2$ is odd, say $k_1-k_2=2t-1$, then $\gamma=\alpha^{i_1-i_2}\beta^{j_1-j_2}\gamma^{2t}$. As $\alpha=\rho^{4m}$, the right hand side is of the form $\rho^x$. On the other hand, $\gamma(A_0)=B_{5m+1}$. Thus $\gamma\neq \alpha^{i_1-i_2}\beta^{j_1-j_2}\gamma^{2t}$. Hence $k_1-k_2$ is even, say $2t$. Thus, we have $\rho^{4m(i_2-i_1)}=\rho^{(4m+4)t+3\frac{m}{2}(j_1-j_2)}, \mbox{ i.e., }$
\begin{equation}\label{family-4-m=12l+6-claim-2-equation-2}
12m \mbox{ divides }4(m+1)t+3(6l+3)(j_1-j_2)+4m(i_1-i_2).
\end{equation}

This implies that $4$ divides $9(2l+1)(j_1-j_2)$, i.e., $4|(j_1-j_2)$. Now as $0\leq j_1-j_2\leq 7$, we have $j_1-j_2=0$ or $4$.\\
\noindent {\it Sub-claim 1a:} $j_1-j_2=0$.\\
If possible, let $j_1-j_2=4$. Then, from Equation \ref{family-4-m=12l+6-claim-2-equation-2}, we have $12m \mbox{ divides }4(m+1)t+6m+4m(i_1-i_2)$ and hence $m|4(m+1)t$, i.e., $m|4t$, as $gcd(m.m+1)=1$. Now, as $0\leq 4t=2(k_1-k_2)\leq 2m-2$, we have $4t=0$ or $m$. However, if $4t=m$, we have $2t=(6l+3)$, an odd number. Thus $4t$ and hence $t=0$. Therefore, from Equation \ref{family-4-m=12l+6-claim-2-equation-2}, we get $12m$ divides $6m+4m(i_1-i_2)$, i.e., $6|4(i_1-i_2)$ which implies $3|(i_1-i_2)$ i.e., $i_1=i_2$. However, this implies that $12m|6m$, a contradiction. Thus Sub-claim 1a is true and $j_1=j_2$. Thus Equation \ref{family-4-m=12l+6-claim-2-equation-2} reduces to 
\begin{equation}\label{family-4-m=12l+6-claim-2-equation-3}
3m \mbox{ divides }(m+1)t+m(i_1-i_2).
\end{equation}
Again since $gcd(m,m+1)=1$,  this implies that $m|t$. However, as $0\leq t\leq \frac{m-1}{2}$, we have $t=0$ and hence $k_1=k_2$. Thus from Equation \ref{family-4-m=12l+6-claim-2-equation-3}, we get $3|(i_1-i_2)$, i.e., $i_1=i_2$. Thus Claim 1 is true and $|H|=24m=2n$. So, as in proof of Claim 2 in Theorem \ref{family-4-m-odd-cayley-theorem}, it suffices to show that $\mathsf{Stab}_H(A_0)=\{\mathsf{id}\}$. Let $\alpha^i\beta^j\gamma^k(A_0)=A_0$. As $\alpha=\rho^{4m},\beta=\rho^{3m/2}$ and $\gamma^2=\rho^{4m+4}$ are powers of $\rho$ and $\gamma(A_0)=B_{5m+1}$, if $k$ is odd, $\alpha^i\beta^j\gamma^k(A_0)=B_x$ for some index $x$. Thus $k$ is even, say $k=2t$. Thus, we have $\alpha^i\beta^j\gamma^k(A_0)=\rho^{4mi+8(m+1)t+\frac{3mj}{2}}=A_0$, i.e., $12m$ divides $4mi+8(m+1)t+\frac{3mj}{2}$, i.e., 
\begin{equation}\label{family-4-m=12l+6-regular-equation-1}
24m\mbox{ divides }8mi+16(m+1)t+3mj
\end{equation}
This implies that $m|16(m+1)t$. As $gcd(m,m+1)=1$, we have $m|16t$. Again, as $m=12l+6=2(6l+3)$ and $6l+3$ is odd, we have $m|2t=k$, i.e., $k=t=0$. Thus Equation \ref{family-4-m=12l+6-regular-equation-1} reduces to $24m$ divides $8mi+3mj$, i.e., $24|(8i+3j)$. However, this implies that $8|j$ and $3|i$, i.e., $i=j=0$. Thus $\mathsf{Stab}_H(A_0)=\{\mathsf{id}\}$ and the theorem holds.\qed

{\theorem \label{family-4-m=12l+10-Cayley-theorem} If $m\equiv 10(mod~12)$, then $R_{12m}(3m+2,9m+1)$ is a Cayley graph.}\\
\pf Let $m=12l+10$. Therefore $8m=96l+80$, i.e., $8m-8=12(8l+6)$. By Proposition \ref{family-4-gamma-results}, we have $\gamma^2=\rho^{12}$. Define $$H=\langle \alpha,\beta,\gamma: \alpha^3=\beta^8=\gamma^{2m}=\mathsf{id};\alpha\beta=\beta\alpha, \alpha\gamma=\gamma\alpha,  \gamma\beta=\beta^{m+1}\gamma, \gamma^m=\alpha^2\beta^4 \rangle.$$
Thus, it is clear that every element of $H$ is of the form $\alpha^i\beta^j\gamma^k$ where $i=0,1,2$; $j=0,1,\ldots,7$ and $k=0,1,\ldots,m-1$. 

\noindent {\it Claim 1:} $H=\{\alpha^i\beta^j\gamma^k:i=0,1,2; j=0,1,\ldots,7; k=0,1,\ldots,m-1\}$.\\
{\it Proof of Claim 1:} If possible, let there exist $i_1,i_2\in \{0,1,2 \}, j_1,j_2 \in \{0,1,\ldots,7\}$ and $k_1,k_2 \in \{0,1,\ldots, m-1\}$, such that $\alpha^{i_1}\beta^{j_1}\gamma^{k_1}=\alpha^{i_2}\beta^{j_2}\gamma^{k_2}$. As $\alpha\beta=\beta\alpha$ and $\alpha\gamma=\gamma\alpha$, we have 
\begin{equation}\label{family-4-m=12l+10-claim-1-equation-1}
\alpha^{i_2-i_1}=\beta^{j_1-j_2}\gamma^{k_1-k_2}.
\end{equation}
If $k_1-k_2$ is odd, say $k_1-k_2=2t-1$, then $\gamma\alpha^{i_2-i_1}=\beta^{j_1-j_2}\gamma^{2t}$. As $\gamma^2=\rho^{12}$, the right hand side is of the form $\rho^x$, i.e., $\gamma\alpha^{i_2-i_1}=\rho^x$. Now $i_2-i_1=0,1$ or $2$. Thus either of $\gamma,\alpha\gamma,\alpha^2\gamma$ is $\rho^x$. But $\gamma(A_0)=B_{9m+5},\alpha\gamma(A_0)=B_{10m+5},~\alpha^2\gamma(B_0)=A_{5m+7}$.
As each of $\gamma,\alpha\gamma,\alpha^2\gamma$ maps some $A_i$ to some $B_j$, none of them is equal to $\rho^x$ and hence a contradiction. So $k_1-k_2$ is even, say $k_1-k_2=2t$. As $\gamma^2=\rho^{12}$ and $\beta=\rho^{3m/2}$, we have $\alpha^{i_2-i_1}=\rho^x$, i.e., $(\rho\sigma)^{4m(i_2-i_1)}=\rho^x$. This implies that $3$ divides $4m(i_2-i_1)$, i.e., $3|m$ or $3|(i_2-i_1)$. As $3$ does not divide $m$, we have $3|(i_2-i_1)$, i.e., $i_1=i_2$. Thus $\rho^{\frac{3m}{2}(j_1-j_2)}=\beta^{j_1-j_2}=\gamma^{k_2-k_1}=(\gamma^2)^t=\rho^{12t}$, i.e., 
\begin{equation}\label{family-4-m=12l+10-claim-1-equation-2}
24m\mbox{ divides }3m(j_1-j_2)-24t
\end{equation}
Thus, we have $m|24t$. As $m=2(6l+5)$, $(6l+5)$ is odd and $3$ does not divide $(6l+5)$, we get $\frac{m}{2}|t$. However, as $0\leq k_2-k_1\leq m-1$, we have $0\leq t \leq \frac{m-1}{2}$. Hence $t=0$ and $k_1=k_2$. Also Equation \ref{family-4-m=12l+10-claim-1-equation-2} reduces to $8|(j_1-j_2)$. Thus $j_1=j_2$. Hence Claim 1 is true and $|H|=24m=2n$. 

So, as in proof of Claim 2 in Theorem \ref{family-4-m-odd-cayley-theorem}, it suffices to show that $\mathsf{Stab}_H(A_0)=\{\mathsf{id}\}$. Let $\alpha^i\beta^j\gamma^k(A_0)=A_0$.\\
\noindent {\it Claim 2:} $k$ is even.\\
{\it Proof of Claim 2:} If possible, let $k$ be odd, say $k=2t+1$. Then, as $\alpha$ commutes with $\beta$ and $\gamma$, we have $\beta^j\gamma^{2t}\gamma\alpha^i(A_0)=A_0$, i.e., $\gamma\alpha^i(A_0)=\beta^{-j}(\gamma^2)^{-t}(A_0)=\rho^x(A_0)=A_x$, as in the proof of Claim 1 of this theorem. Now, $i=0,1$ or $2$ and as $\gamma(A_0)=B_{9m+5}$ and $\alpha\gamma(A_0)=B_{10m+5}$, we have $i=2$. This implies $\alpha^2\beta^j\gamma^{2t+1}(A_0)=A_0$, i.e., $\beta^j(\gamma^2)^t\gamma(A_0)=\alpha(A_0)=A_{7m}$, i.e., $$A_{7m}=\beta^j(\gamma^2)^t\gamma(A_0)=\beta^j(\gamma^2)^t(B_{9m+5})=\rho^x(B_{9m+5})=B_{9m+x+5}, \mbox{ a contradiction.}$$ Hence the claim is true and let $k=2t$. Therefore, $$\beta^j(\gamma^2)^t(A_0)=\alpha^{-i}(A_0).$$ As left side of the above equation is $\rho^x(A_0)$ and $\alpha^2(A_0)=B_{2m-1}$, we conclude that $i=0$ or $2$. If $i=2$, then we have $\alpha^2\beta^j(\gamma^2)^t(A_0)=A_0$. Again as $\alpha$ commutes with $\beta$ and $\gamma$, we have $$A_0=\beta^j\gamma^{2t}\alpha^2(A_0)=\beta^j\gamma^{2t}(B_{2m-1})=\rho^x(B_{2m-1})=B_{2m+x-1}, \mbox{ a contradiction}.$$
Therefore, $i=0$ and hence we have $\beta^j(\gamma^2)^t(A_0)=A_0$, i.e., $$\rho^{12t+3j\frac{m}{2}}(A_0)=A_0,\mbox{ i.e., }12m \mbox{ divides }12t+3j\frac{m}{2}=12t+3j(6l+5)$$
Thus $12$ divides $3j(6l+5)$, i.e., $4|j(6l+5)$. However as $6l+5$ is odd and $j\in \{0,1,\ldots,7\}$, we have $j=0$ or $4$. If $j=4$, we have $12m$ divides $12t+12(6l+5)$, i.e., $m=12l+10=2(6l+5)$ divides $t+(6l+5)$ and hence $(6l+5)$ divides $t$. However as $0\leq k \leq m-1$, we have $0\leq t \leq \frac{m-1}{2}<6l+5$. Thus the only possible value of $t$ is $0$ and hence $k=0$. Therefore, we have $\beta^j(A_0)=A_0$, i.e., $\rho^{3(6l+5)j}(A_0)=A_0$. This implies that $12m=24(6l+5)$ divides $3(6l+5)j$, i.e., $8|j$ and hence $j=0$. 

Thus we have $\mathsf{Stab}_H(A_0)=\{\mathsf{id}\}$ and the theorem holds.\qed

\section{Family-5 [$R_{2m}(2b,r)$: $b^2\equiv \pm 1$(mod $m$) and $r \in \{1,m-1\}$ is odd]}
{\theorem If $b^2\equiv \pm 1$(mod $m$) and $r \in \{1,m-1\}$ is odd, then $R_{2m}(2b,r)$ is a Cayley graph.}\\
\pf If $r=1$, then it is clear that the conditions of being in {\bf Family-1} are satisfied, (i.e., $r^2\equiv 1(mod~n)$ and $ra\equiv a(mod~n)$) and hence, by Theorem \ref{r^2=1}, $R_{2m}(2b,r)$ is a Cayley graph. So we are left with the case when $n=2m$, $a=2b$, $b^2\equiv \pm 1(mod~m)$, $r=m-1$ and $m$ is even. Observe that, in this case, $$r^2=(m-1)^2=m^2-2m+1\equiv 1(mod~2m)\equiv 1(mod~n)~[\mbox{since, }m \mbox{ is even}].$$
Also, as $m$ divides $bm$ i.e., $m|b(r+1)$, we have $br\equiv -b(mod~m)$,  i.e., $2br\equiv -2b(mod~2m)$, i.e., $ra\equiv -a(mod~n)$. Thus, in this case, $r^2\equiv 1(mod~n)$ and $ra\equiv -a(mod~n)$ holds. Hence, by Theorem \ref{r^2=1}, $R_{2m}(2b,r)$ is a Cayley graph.\qed 

{\remark The above theorem shows that {\bf Family-5} is a subfamily of {\bf Family-1}. However, they were shown as different families in Theorem 3.10 in \cite{vt-rw}. }

Combining the analysis of the rose window graphs in {\bf Families: 1--5}, we have Theorem \ref{cayley-iff}.

\section{Appendix}

{\lemma \label{family-4-calc} Let $G=\mathsf{Aut}(R_{12m}(3m+2,9m+1))$ , where $m\equiv 0(mod~4)$. Then $|G|=96m$.}\\
\pf Since, $R_{12m}(3m+2,9m+1)$ is vertex-transitive and its order is $24m$ and $\mathsf{Stab}_G(A_0)$ contains $\mathsf{id}, \mu,\sigma,\mu\sigma$, therefore, by orbit-stabilizer theorem, we have $|G|\geq 4\times 24m=96m$. Thus, it is enough to show that $|G|\leq 96m$. We also  know that 
$$G=\langle \rho,\mu,\sigma: \rho^n=\mu^2=\sigma^2=\mathsf{id}; \mu\rho\mu=\rho^{-1}, \sigma\rho^3\sigma=\rho^3,\sigma\mu=\mu\sigma,~~~~~~~~~~~~~~~~~~$$ 
$$~~~~~~~(\rho\sigma)^3=(\sigma\rho)^3=\rho^{3(m+1)}, (\rho\sigma\rho)^3=\rho^{9m+6} \rangle, \mbox{ where }n=12m.$$
Consider the sets $X=\{\rho^i\sigma\rho^j\mu^k:i\in \{0,1,2\ldots,n-1\},j \in \{0,1,2\},k  \in \{0,1\}\}$ and $Y=\{\rho^i\mu^k: i\in \{0,1,2\ldots,n-1\},k  \in \{0,1\}\}$. We claim that all elements are either in $X$ or in $Y$. It is clear that elements in $G$ which does not involve $\sigma$ are in $Y$, due to the relations $\rho^n=\mu^2=\mathsf{id}$ and $\mu\rho\mu=\rho^{-1}$. Again, as $\sigma\mu=\mu\sigma$ and $\mu\rho=\rho^{-1}\mu$, any element in $G$ can be expressed in the form where $\mu$ occurs in the extreme right of the expression. Thus it is enough to show that elements in $G$ which involve only $\rho$ and $\sigma$ are of the form $\rho^i\sigma\rho^j$ where $i\in \{0,1,2\ldots,n-1\}$ and $j \in \{0,1,2\}$. Again, as $\sigma\rho^3=\rho^3\sigma$, it is clear that the power of $\rho$ lying on the right of $\sigma$ can be made $0,1$ or $2$. Finally, we deal with elements $\sigma\rho\sigma$ and $\sigma\rho^2\sigma$. 

As $(\rho\sigma\rho)^3=\rho^{9m+6}$, we have $\sigma\rho^2\sigma\rho^2\sigma=\rho^{9m+4}$, i.e., $$\sigma\rho^2\sigma=\rho^{9m+4}\sigma\rho^{-2}=\rho^{9m+4}\sigma\rho^{12m-2}=\rho^{9m+4+12m-3}\sigma\rho=\rho^{9m+1}\sigma\rho \in X.$$

As $(\rho\sigma)^3=\rho^{3(m+1)}$, we have $(\sigma\rho\sigma\rho\sigma)=\rho^{3m+2}$, i.e., $$\sigma\rho\sigma=\rho^{3m+2}\sigma\rho^{-1}=\rho^{3m+2}\sigma\rho^{12m-1}=\rho^{3m+2+12m-3}\sigma\rho^2=\rho^{3m-1}\sigma\rho^{2}\in X.$$

Similarly, any other element of $G$ involving $\rho$ and $\sigma$ can be expressed in the form of elements in $X$. Thus $G=X \cup Y$ and hence $$|G|=|X \cup Y|\leq |X|+|Y|\leq (n\times 3 \times 2)+(n \times 2)=6n+2n=8n=96m.$$\qed

\noindent {\bf Proof of Proposition \ref{family-4-gamma-results}	}: 
\begin{enumerate}
	\item For $m=12l+2$, we have $8m=96l+16$, i.e., $8m-4=12(8l+1)$.  $$\gamma^2=(\rho^{8m}\sigma\rho^2\omega)(\rho^{8m}\sigma\rho^2\omega)=\rho^{8m}\rho^{8m-4}\sigma\rho^2\omega\rho^{4}\sigma\rho^2\omega~~~~(\mbox{as }\rho^{12} \mbox{ commutes with }\sigma \mbox{ and } \omega)$$
	$$~~~=\rho^{4m-4}\sigma\rho^2(\omega\rho^{3})\rho\sigma\rho^2\omega=\rho^{4m-4}\sigma\rho^2(\rho^{3(m+1)}\omega)\rho\sigma\rho^2\omega ~~~~~~~~~~~~~~~~(\mbox{as }\omega\rho^{3l}=\rho^{3l(m+1)})$$
	$$=\rho^{7m-1}\sigma\rho^2(\omega\rho)\sigma\rho^2\omega=\rho^{7m-1}\sigma\rho^2(\sigma\rho\omega)\sigma\rho^2\omega=\rho^{7m-1}\sigma\rho^2\sigma\rho\sigma\omega\rho^2\omega~~~~~~~~~~~~~~~~~~~~$$
	$$=\rho^{7m-1}\sigma\rho^2\sigma\rho\sigma(\omega\rho\omega)^2=\rho^{7m-1}\sigma\rho^2\sigma\rho\sigma(\sigma\rho)^2=\rho^{7m-1}\sigma\rho^2\sigma\rho\sigma(\sigma\rho)(\sigma\rho)~~~~~~~~~~~~~~~$$
	$$=\rho^{7m-1}\sigma\rho^2\sigma\rho^2\sigma\rho=\rho^{7m-2}\rho\sigma\rho^2\sigma\rho^2\sigma\rho=\rho^{7m-2}(\rho\sigma\rho)(\rho\sigma\rho)(\rho\sigma\rho)=\rho^{7m-2}(\rho\sigma\rho)^3~~~~$$ 
	$$=\rho^{7m-2}\rho^{9m+6}=\rho^{16m+4}=\rho^{4m+4}~~~~~~~~~~~~~~~~~~~~~~~~~~~~~~~~~~~~~~~~~~~~~~~~~~~~~~~~~~~~~~~~$$
	For $m=12l+6$, we have $8m=96l+48=12(8l+4)$.  $$\gamma^2=(\rho^{8m}\sigma\rho^2\omega)(\rho^{8m}\sigma\rho^2\omega)=\rho^{16m}\sigma\rho^2\omega\sigma\rho^2\omega~~~~(\mbox{as }\rho^{12} \mbox{ commutes with }\sigma \mbox{ and } \omega)$$
	$$~~~=\rho^{4m}\sigma\rho^2\sigma\omega\rho^2\omega=\rho^{4m}(\sigma\rho\sigma)^2(\omega\rho\omega)^2=\rho^{4m}(\sigma\rho\sigma)^2(\sigma\rho)^2 ~~~~~~~~~(\mbox{as }\omega\rho=\sigma\rho\omega)$$
	$$=\rho^{4m}\sigma\rho^3\sigma\rho=\rho^{4m+4}.~~~~~~~~~~~~~~~~~~~~~~~~~~~~~~~~~~~~~~~~~~~~~~~~~~~~~~~~~~~~~~~~~~~$$
	Similarly, for $m=12l+10$, it can be proved that $\gamma^2=\rho^{12}$.
	\item The values of $\gamma^m$ can be found by raising $\gamma^2$ to the power $m/2$, and hence can be checked to have the respective forms. 
\end{enumerate}

\noindent {\bf Checking whether a rose window graph is Cayley using SageMath}

The following is the code to check whether a rose window graph is Cayley. The code is given for $R_{36}(11,28)$, which was claimed to be Cayley in Theorem \ref{family-4-m=3-theorem}. Readers can also edit the values of $n,a,r$ to check for other rose window graphs. The output will be TRUE, if the graph is Cayley, else it will be FALSE.
\begin{verbatim}
          n=36
          a=11
          r=28
          A = list(var('A_%d' % i) for i in range(n))
          B = list(var('B_%d' % i) for i in range(n))
          V=A+B
          E=[]
          G=Graph()
          G.add_vertices(V)
          for i in range(n):
              E.append((A[i],A[mod(i+1,n)]))
              E.append((A[i],B[i]))
              E.append((B[i],A[mod(i+a,n)]))
              E.append((B[i],B[mod(i+r,n)]))
          G.add_edges(E)
          G.is_cayley()
\end{verbatim}

\end{document}